\documentclass{article}
\usepackage{eurosym}
\usepackage{amssymb}
\usepackage{amsfonts}
\usepackage{graphicx}
\usepackage{amsmath}

\newtheorem{theorem}{Theorem}

\newtheorem{corollary}[theorem]{Corollary}

\newtheorem{definition}[theorem]{Definition}

\newtheorem{lemma}[theorem]{Lemma}

\newtheorem{proposition}[theorem]{Proposition}
\newtheorem{remark}[theorem]{Remark}

\newcommand{\LL}{{\mathbb L}}

\renewcommand{\SS}{{\mathbb S}}
\renewcommand{\P}{{\mathbb P}}

\DeclareMathOperator{\Mon}{Mon}

\begin{document}

\begin{center}
{\LARGE Law of the iterated logarithm for the periodogram} \bigskip

Christophe Cuny$^{a}$, Florence Merlev\`{e}de$^{b}$ and Magda Peligrad$^{c}$%
\footnote{\textit{\ }Supported in part by a Charles Phelps Taft Memorial
Fund grant, the NSA\ grant H98230-11-1-0135 and the NSF\ grant DMS-1208237.}
\end{center}

$^{a}$ Laboratoire MAS, Ecole Centrale de Paris, Grande Voie des Vignes,
92295 Chatenay-Malabry cedex, FRANCE. E-mail: christophe.cuny@ecp.fr

$^{b}$ Universit\'e Paris-Est, LAMA (UMR 8050), UPEMLV, CNRS, UPEC, F-77454
Marne-La-Vall\'ee, FRANCE. E-mail florence.merlevede@univ-mlv.fr

$^{c}$ Department of Mathematical Sciences, University of Cincinnati, PO Box
210025, Cincinnati, Oh 45221-0025, USA. E-mail: peligrm@ucmail.uc.edu

\medskip



\bigskip

\textbf{Abstract }

\bigskip

We consider the almost sure asymptotic behavior of the periodogram of
stationary and ergodic sequences. Under mild conditions we establish that the limsup of the
periodogram properly normalized identifies almost surely the spectral
density function associated with the stationary process. Results for a
specified frequency are also given. Our results also lead to the law of the
iterated logarithm for the real and imaginary part of the discrete Fourier
transform. The proofs rely on martingale approximations combined with
results from harmonic analysis and technics from ergodic theory. Several
applications to linear processes and their functionals, iterated random
functions, mixing structures and Markov chains are also presented.

\section{Introduction}

The periodogram, introduced as a tool by Schuster in 1898, plays an
essential role in the estimation of the spectral density of a stationary
time series $(X_{j})_{j\in \mathbb{Z}}$ of centered random variables with
finite second moment. The finite Fourier transform is defined as%
\begin{equation}
S_{n}(t)=\sum_{k=1}^{n}\mathrm{e}^{ikt}X_{k}\,,  \label{Four}
\end{equation}%
where $i=\sqrt{-1}$ is the imaginary unit, and the periodogram as 
\begin{equation}
I_{n}(t)=\frac{1}{2\pi n}\big |S_{n}(t)\big |^{2}\ \,\ t\in \lbrack 0,2\pi
]\,.  \label{Per}
\end{equation}

It is well-known since Wiener and Wintner \cite{Win} that for any stationary
sequence $(X_{j})_{j\in \mathbb{Z}}$ in ${\mathbb{L}}^{1}$ (namely ${\mathbb{%
E}}|X_{0}|<\infty $) there is a set $\Omega ^{\prime }$ of probability one
such that for any $t\in \lbrack 0,2\pi ]$ and any $\omega \in \Omega
^{\prime },$ $S_{n}(t)/n$ converges. To provide the speed of this convergence many authors
(see Peligrad and Wu \cite{PeWu} and the references therein) established a
central limit theorem for the real and imaginary parts of ${\ S_{n}(}${$t$}${%
)}/\sqrt{n}$ under various assumptions$.$ Recently, Peligrad and Wu \cite%
{PeWu} showed that, under a very mild regularity condition and finite second
moment, $\frac{1}{\sqrt{n}}[{\mathcal{R}e}(S_{n}(t)),\,{\mathcal{I}m}%
(S_{n}(t))]$ are asymptotically independent normal random variables with
mean $0$ and variance $\pi f(t)$ for almost all $t$ (here $f$ is the
spectral density of $(X_{j})_{j\in \mathbb{Z}}$). The central limit theorem
implies that $I_{n}(t)/\log \log n$ converges to $0$ in probability. An
interesting and natural problem, that apparently has never been studied in
depth before, is the law of the iterated logarithm, namely to identify in
the almost sure sense, $\limsup_{n\rightarrow \infty }I_{n}(t)/\log \log n$
for almost all $t$, or for a $t$ fixed. In this paper, we study both these
problems. We provide mild sufficient conditions on the stationary sequence that are sufficient to have $%
\limsup_{n\rightarrow \infty }I_{n}(t)/\log \log n=f(t)$ almost surely.
These results shed additional light on the importance of the periodogram in
approximating the spectral density $f(t)$ of a stationary process. The
techniques are based on martingale approximation, rooted in Gordin \cite{G1}
and Rootz\'{e}n \cite{Rootzen} and developed by Gordin and Lifshitz \cite{GL}
and Woodroofe \cite{Wood}, combined with tools from ergodic theory and
harmonic analysis. Various applications are presented to linear processes
and their functionals, iterated random functions, mixing structures and
Markov chains.

We would like to point out that our results are formulated under the
assumption that the underlying stationary sequence is assumed to be adapted
to an increasing (stationary) filtration. Results in the non adapted case
could also be obtained. We shall also assume that our stationary sequence is
constructed via a measure-preserving transformation that is invertible.
Since our proofs are based on martingale approximation, we could also obtain
similar results when the measure-preserving transformation is assumed to be
non invertible. In this situation, the conditions should be expressed with
the help of the Perron-Frobenius operator associated to the transformation
(see for instance \cite{CM}). In our paper we shall not pursue these last
two cases.

Our paper is organized as follows. Section \ref{sectres} contains the
presentation of the results. Section \ref{sectproofs} is devoted to the
proofs. Applications are presented in Section \ref{sectappl}.

\section{Main results}

\label{sectres}

Let $(\Omega ,{\mathcal{A}},{\mathbb{P}})$ be a probability space. Assume,
without loss of generality, that ${\mathcal{A}}$ is a countably generated $%
\sigma $-field, and let $\theta :\Omega \rightarrow \Omega $ be a bijective
bi-measurable transformation preserving ${\mathbb{P}}$. Let ${\mathcal{F}}%
_{0}$ be a $\sigma $-algebra such that ${\mathcal{F}}_{0}\subseteq \theta
^{-1}({\mathcal{F}}_{0})$. Let $({\mathcal{F}}_{n})_{n\in \mathbb{Z}}$ be
the non-decreasing filtration given by ${\mathcal{F}}_{n}=\theta ^{-n}({%
\mathcal{F}}_{0})$, and let ${\mathcal{F}}_{-\infty }=\bigcap_{k\in {\mathbb{%
Z}}}{\mathcal{F}}_{k}$. All along the paper $X_{0}$ is a centered real random variable in ${\mathbb{L}}^{2}$ which is ${%
\mathcal{F}}_{0}$-measurable. We then define a stationary sequence $%
(X_{n},n\in {\mathbb{Z}})$ by $X_{n}=X_{0}\circ \theta ^{n}$. We denote ${%
\mathbb{E}}_{k}(\cdot )={\mathbb{E}}(\cdot |\mathcal{F}_{k})$ and $%
P_{k}(\cdot )={\mathbb{E}}_{k}(\cdot )-{\mathbb{E}}_{k-1}(\cdot )$.
Throughout the paper, we say that a complex number $z$ is an eigenvalue of $%
\theta $ if there exists $h\neq 0$ in ${\mathbb{L}}^{2}({\mathbb{P}})$ such
that $h\circ \theta =zh$ almost everywhere. We say that $A\in {\mathcal{A}}$
is invariant if $\theta ^{-1}(A)=A$. If for any invariant set $A$, ${\mathbb{%
P}}(A)=0$ or $1$, we say that $\theta $ is ergodic with respect to ${\mathbb{%
P}}$, or equivalently that the stationary sequence is ergodic.

Relevant to our results is the notion of spectral distribution function
induced by the covariances. By Herglotz's Theorem (see e.g. Brockwell and
Davis \cite{BD}), there exists a non-decreasing function $G$ (the spectral
distribution function) on $[0,2\pi ]$ such that, for all $j\in \mathbb{Z}$, 
\begin{equation*}
\mathrm{cov}(X_{0},X_{j})=\int_{0}^{2\pi }\exp (ij\theta )dG(\theta ),\quad
j\in \mathbb{Z}\,.
\end{equation*}%
If $G$ is absolutely continuous with respect to the normalized Lebesgue
measure $\lambda $ on $[0,2\pi ]$, then the Radon-Nikodym derivative $f$ of $%
G$ with respect to the Lebesgue measure is called the spectral density and
we have 
\begin{equation*}
\mathrm{cov}(X_{0},X_{j})=\int_{0}^{2\pi }\exp (ij\theta )f(\theta )d\theta
,\quad j\in \mathbb{Z}\,.
\end{equation*}

Our first theorem points out a projective condition which assures the law of
the iterated logarithm for almost all frequencies.

All along the paper, denote 
\begin{equation*}
Y_{k}(t)=\big (\cos (kt)X_{k},\sin (kt)X_{k}\big )^{\prime }\, ,
\end{equation*}
where $u^{\prime }$ stands for the transposed vector of $u$.

\begin{theorem}
\label{thmlilmoment2}Assume that $\theta $ is ergodic and that 
\begin{equation}
\sum_{k \geq 2}\frac{(\log k )}{k }\Vert {\mathbb{E}}_{0}(X_{k })\Vert
_{2}^{2}<\infty \,.  \label{condFM2}
\end{equation}%
Then the spectral density, say $f$, of $(X_{k},k\in {\mathbb{Z}})$ exists
and for almost all $t\in \lbrack 0,2\pi )$, the sequence $\big \{%
\sum_{k=1}^{n}Y_{k}(t)/\sqrt{2n\log \log n},n\geq 3\big \}$ is $%
\mbox{${\mathbb P}$-a.s.}$ bounded and has the ball $\{x\in {\mathbb{R}}%
^{2}\,:\,x^{\prime }x\leq \pi f(t)\}$ as its set of limit points. In
particular, for almost all $t\in \lbrack 0,2\pi )$, the following law of the
iterated logarithm holds 
\begin{equation*}
\limsup_{n\rightarrow \infty }\frac{I_{n}(t)}{\log \log n}=f(t)\text{ }\ %
\mbox{${\mathbb P}$-a.s.}
\end{equation*}
\end{theorem}

Note that condition (\ref{condFM2}) is satisfied by martingale differences.
It is a very mild condition involving only a logarithmic rate of convergence
to $0$ of $\Vert {\mathbb{E}}_{0}(X_{k })\Vert _{2}.$

\bigskip

If we assume a more restrictive moment condition, \eqref{condFM2} can be
weakened. Define the function $L(x)=\log (\mathrm{e}+|x|)$.

\begin{theorem}
\label{thmlilmoment2delta}Assume that $\theta $ is ergodic. Assume in
addition that 
\begin{equation*}
{\mathbb{E}}\Big (\frac{X_{0}^{2}L(X_0)}{L(L(X_{0}))}\Big )<\infty \, ,
\end{equation*}
and that 
\begin{equation}
\sum_{k \geq 3}\frac{\Vert {\mathbb{E}}_{0}(X_{k })\Vert _{2}^{2}}{k (\log
\log k )}<\infty \,.  \label{condFM}
\end{equation}%
Then the conclusions of Theorem \ref{thmlilmoment2} hold.
\end{theorem}

Note that condition \eqref{condFM2}, as well as condition \eqref{condFM},
implies the following regularity condition: 
\begin{equation}
{\mathbb{E}}(X_{0}|{\mathcal{F}}_{-\infty })=0\ \quad 
\mbox{${\mathbb
P}$-a.s.}  \label{regular}
\end{equation}%
We point out that this regularity condition implies that the process $%
(X_{k})_{k\in {\mathbb{Z}}}$ is purely non deterministic. Hence by a result
of Szeg\"{o} (see for instance \cite[Theorem 3]{Bingham}) if \eqref{regular}
holds, the spectral density $f$ of $(X_{k})_{k\in {\mathbb{Z}}}$ exists and
if $X_{0}$ is non degenerate, 
\begin{equation*}
\int_{0}^{2\pi }\log f(t)~dt>-\infty \,;
\end{equation*}%
in particular, $f$ cannot vanish on a set of positive measure. We mention
also that under \eqref{regular}, Peligrad and Wu \cite{PeWu} established
that 
\begin{equation}
\lim_{n\rightarrow \infty }\frac{{\mathbb{E}}|S_{n}(t)|^{2}}{n}=2\pi f(t)\ \ 
\text{for almost all $t\in \lbrack 0,2\pi )$}\,  \label{defgt}
\end{equation}%
(see their Lemma 4.2).

\smallskip

Both theorems above
hold for almost all frequencies. It is possible that on a set of measure $0$
the behavior be quite different. This fact is suggested by a result of
Rosenblatt \cite{Ro81} who established, on a set of measure $0$, non-normal
attraction for the Fourier transform under a different normalization than $%
\sqrt{n}$.

\bigskip

We give next conditions imposed to the stationary sequence which help to
identify the frequencies for which the LIL\ holds. As we shall see, the next
result is well adapted for linear processes generated by iid (independent
identically distributed) sequences.

\begin{theorem}
\label{thm3}Assume that \eqref{regular} holds and that%
\begin{equation}
\sum_{n\geq 0}\Vert P_{0}(X_{n})-P_{0}(X_{n+1})\Vert _{2}<\infty \,.
\label{condproj}
\end{equation}%
Then the spectral density $f(t)$ of $(X_{k}, k \in {\mathbb{Z}})$ is
continuous on $(0, 2 \pi)$, and the convergence \eqref{defgt} holds for all $%
t\in (0,2\pi )$. Moreover if $\theta $ is ergodic, the conclusions of
Theorem \ref{thmlilmoment2} hold for all $t\in (0,2\pi )\backslash \{\pi \}$
such that $\mathrm{e}^{-2it}$ is not an eigenvalue of $\theta $.
\end{theorem}

\begin{remark}
The conditions of this theorem do not imply that the spectral density is
continuous at $0$. This is easy to see by considering the time series $%
X_{k}=\sum_{j\geq 0}j^{-3/4}\varepsilon _{k-j}$ where $(\varepsilon
_{k})_{k\in {\mathbb{Z}}}$ is a sequence of iid centered real random
variables in ${\mathbb{L}}^{2}.$ For this case all the conditions of Theorem %
\ref{thm3} are satisfied (see Section \ref{sectionLP}) and \textrm{var}$%
(\sum_{k=1}^{n}X_{k})/n$ converges to $\infty $. This shows that the
spectral density is not continuous at $0$ since otherwise we would have $%
\mathrm{Var}(\sum_{k=1}^{n}X_{k})/n\rightarrow 2\pi f(0),$ which is not the
case.
\end{remark}

\medskip

We would like to mention that Condition (\ref{condproj}) above was used by
Wu \cite{Wu05} in the context of the CLT. We also infer from our proof and
Remark \ref{remark0pi}, that if \eqref{regular} and \eqref{condproj} hold,
and $\theta $ is ergodic, then $\limsup_{n\rightarrow \infty }I_{n}(\pi
)/\log \log n=2f(\pi )$ ${\mathbb{P}}$-a.s. Moreover, it follows from a recent result of Cuny \cite{CunyLIL} that if 
$\theta $ is ergodic, \eqref{regular} holds and \eqref{condproj} is
reinforced to $\sum_{n\geq 0}\Vert P_{0}(X_{n})\Vert _{2}<\infty $, then $%
\limsup_{n\rightarrow \infty }I_{n}(0)/\log \log n=2f(0)$ ${\mathbb{P}}$-a.s.

\bigskip

We say that $\theta $ is \textit{weakly-mixing}, if for all $A,B\in {%
\mathcal{A}}$, 
\begin{equation*}
\lim_{n\rightarrow \infty }\frac{1}{n}\sum_{k=0}^{n-1}|{\mathbb{P}}(\theta
^{-k}A\cap B)-{\mathbb{P}}(A){\mathbb{P}}(B)|=0\,.
\end{equation*}%
It is well-known (see e.g. \cite[Theorem 6.1]{Petersen}) that saying that $%
\theta $ is weakly mixing is equivalent to saying that $\theta $ is ergodic
and its only eigenvalue is $1$. Let us also mention that when ${\mathcal{F}}%
_{-\infty }$ is trivial then $\theta $ is weakly mixing (see section 2 of 
\cite{Petersen}).

As an immediate corollary to Theorem \ref{thm3} we obtain the following LIL
for all frequencies.

\begin{corollary}
\label{cor5} Assume that $\theta $ is weakly mixing and that \eqref{regular}
and (\ref{condproj}) hold. Then the conclusion of Theorem \ref{thm3} holds
for all $t\in (0,2\pi )\backslash \{\pi \}$.
\end{corollary}

Next theorem involves a projective condition in the spirit of Rootz\'{e}n 
\cite{Rootzen}. It is very useful in order to treat several classes of
Markov chains including reversible Markov chains.

\begin{theorem}
\label{thm4}Assume that $\theta $ is ergodic. Let $t\in (0,2\pi )\backslash
\{\pi \}$ be such that $\mathrm{e}^{-2it}$ is not an eigenvalue of $\theta $%
. Assume in addition that 
\begin{equation}
\sup_{n}\Vert \mathbb{E}_{0}(S_{n}(t))\Vert _{2}<\infty \,.  \label{poisson}
\end{equation}%
Then 
\begin{equation}
\lim_{n\rightarrow \infty }\frac{{\mathbb{E}}|S_{n}(t)|^{2}}{n}=\sigma
_{t}^{2}\,\,\text{ (say)}  \label{defsigmat}
\end{equation}%
and $\big \{\sum_{k=1}^{n}Y_{k}(t)/\sqrt{2n\log \log n},n\geq 3\big \}$ is $%
\mbox{${\mathbb P}$-a.s.}$ bounded and has the ball $\{x\in {\mathbb{R}}%
^{2}\,:\,x^{\prime }x\leq \sigma _{t}^{2}/2\}$ as its set of limit points.
In particular, 
\begin{equation*}
\limsup_{n\rightarrow \infty }\frac{I_{n}(t)}{\log \log n}=\frac{\sigma
_{t}^{2}}{2\pi }\text{ }\ \mbox{${\mathbb P}$-a.s.}
\end{equation*}
\end{theorem}

\begin{remark}
Note that in Theorem \ref{thm4} we do not require the sequence to be
regular, i.e. it may happen that $\mathbb{E}(X_{0}|\mathcal{F}_{-n})$ does
not converge to $0$ in $\mathbb{L}^{2}$. The spectral density might not
exist.
\end{remark}

\section{Proofs}

\label{sectproofs}

\noindent \textbf{Proof of Theorem \ref{thmlilmoment2}.} The proof is based
on martingale approximation. By Lemma 4.1 in Peligrad and Wu \cite{PeWu},
since $X_{0}\in {\mathbb{L}}^{2}$ and \eqref{regular} holds under %
\eqref{condFM2}, we know that for almost all $t\in [0,2\pi )$, the following
limit exists in ${\mathbb{L}}^{2}({\mathbb{P}})$ and $%
\mbox{${\mathbb
P}$-a.s.}$ 
\begin{equation}
D_{0}(t)=\lim_{n\rightarrow \infty }\sum_{k=0}^{n}\mathrm{e}%
^{ikt}P_{0}(X_{k})\,.  \label{defD0alpha}
\end{equation}%
Hence setting for all $\ell \in {\mathbb{Z}}$, 
\begin{equation}
D_{\ell }(t)=\mathrm{e}^{i\ell t}D_{0}(t)\circ \theta ^{\ell }\,,
\label{defDlalpha}
\end{equation}%
we get that for almost all $t\in (0,2\pi )$, $(D_{\ell }(t))_{\ell \in {%
\mathbb{Z}}}$ forms a sequence of martingale differences in ${\mathbb{L}}%
^{2}({\mathbb{P}})$ with respect to $({\mathcal{F}}_{\ell })_{\ell \in {%
\mathbb{Z}}}$. As we shall see, the conclusion of the theorem will then
follow from Propositions \ref{compactlilmart} and \ref{lmaasmartapprox}
below.

\begin{proposition}
\label{compactlilmart}Assume that $\theta $ is ergodic. Let $t\in (0,2\pi
)\backslash \{\pi \}$ and assume that $\mathrm{e}^{-2it}$ is not an
eigenvalue of $\theta .$ Let $D$ be a square integrable complex-valued
random variable adapted to $\mathcal{F}_{0}$ and such that ${\mathbb{E}}%
_{-1}(D)=0$ a.s. For any $k\in {\mathbb{Z}}$, let $d_{k}(t)=\big ({\mathcal{R%
}e}(\mathrm{e}^{ikt}D\circ \theta ^{k}),{\mathcal{I}m}(\mathrm{e}%
^{ikt}D\circ \theta ^{k})\big )^{\prime }$. Then the sequence $\big \{%
\sum_{k=1}^{n}d_{k}(t)/\sqrt{2n\log \log n},n\geq 3\big \}$ is $%
\mbox{${\mathbb P}$-a.s.}$ bounded and has the ball $\{x\in {\mathbb{R}}%
^{2}\,:\,x^{\prime }\cdot x\leq {\mathbb{E}}(|D|^{2})/2\}$ as its set of
limit points.
\end{proposition}

\begin{remark}
\label{remcountably} Since we assume ${\mathcal{A}}$ to be countably
generated, then ${\mathbb{L}}^{2}(\Omega ,{\mathcal{A}},\mathbb{P})$ is
separable and (see Lemma \ref{LLNFourier}) $\theta $ can admit at most
countably many eigenvalues. Hence, Proposition \ref{compactlilmart} applies
to almost all $t\in \lbrack 0,2\pi )$.
\end{remark}

\begin{remark}
\label{remark0pi} Let $t=0$ or $t=\pi $, and assume that $\theta $ is
ergodic. Then if $D$ is a square integrable real-valued random variable
adapted to $\mathcal{F}_{0}$ and such that ${\mathbb{E}}_{-1}(D)=0$ a.s.,
the following result holds: $\limsup_{n\rightarrow \infty
}|\sum_{k=1}^{n}\cos (kt)D\circ \theta ^{k}|^{2}/(2n\log \log n)={\mathbb{E}}%
(D^{2})$ a.s. For $t=0$, it is the usual law of the iterated logarithm for
stationary ergodic martingale differences sequences. For $t=\pi $, it follows from a direct application of \cite[%
Theorem 1]{HS}.
\end{remark}

\begin{proposition}
\label{lmaasmartapprox} Assume that condition \eqref{condFM2} holds. Then,
for almost all $t\in [0,2\pi )$, 
\begin{equation}
\frac{|S_{n}(t)-M_{n}(t)|}{\sqrt{n\log \log n}}\rightarrow 0\ \ %
\mbox{${\mathbb P}$-a.s.}  \label{approxmartlil}
\end{equation}%
where $M_{n}(t)=\sum_{k=1}^{n}D_{k}(t)$ and $D_{k}(t)$ is defined by %
\eqref{defDlalpha}.
\end{proposition}

To end the proof of Theorem \ref{thmlilmoment2}, we proceed as follows. By
Proposition \ref{lmaasmartapprox}, it suffices to prove that the conclusion
of Theorem \ref{thmlilmoment2} holds replacing $Y_{k}(t)$ with $d_{k}(t)=%
\big ({\mathcal{R}e}(D_{k}(t)),{\mathcal{I}m}(D_{k}(t))\big )^{\prime }$.
With this aim, it suffices to apply Proposition \ref{compactlilmart}
together with Remark \ref{remcountably} and to notice the following fact:
according to Lemma 4.2 in Peligrad and Wu \cite{PeWu}, for almost all $t\in
\lbrack 0,2\pi )$, 
\begin{equation*}
\frac{{\mathbb{E}}(|D_{0}(t)|^{2})}{2}=\pi f(t)\,.
\end{equation*}

It remains to prove the above propositions.

\bigskip

It is convenient to work on the product space. Let $(\widetilde{\Omega },%
\widetilde{{\mathcal{F}}},\widetilde{\mathbb{P}})=([0,2\pi ]\times \Omega ,{%
\mathcal{B}}\otimes {\mathcal{A}},\lambda \otimes {\mathbb{P}})$ where $%
\lambda $ is the normalized Lebesgue measure on $[0,2\pi ]$, and ${\mathcal{B%
}}$ be the Borel $\sigma $-algebra on $[0,2\pi ]$. Let $t\in \lbrack 0,2\pi
) $ be a real number, fixed for the moment. Clearly, the transformation $%
\widetilde{\theta}=\widetilde{\theta}_{t}$ (we omit the dependence with
respect to $t $ when $t$ is fixed) given by 
\begin{equation}  \label{thetatilde}
\widetilde{\theta }~:~(u,\omega )\mapsto (u+t\text{ modulo }2\pi ,\theta
(\omega )) \, ,
\end{equation}%
is invertible, bi-measurable and preserves $\widetilde{\mathbb{P}}$.

Consider also the filtration $(\widetilde {\mathcal{F}}_n)_{n\in \mathbb{Z}}$
given by $\widetilde {\mathcal{F}}_n:= {\mathcal{B}}\otimes {\mathcal{F}}_n$.

Define a random variable $\widetilde{X}_{0}$ on $\widetilde{\Omega }$ by $%
\widetilde{X}_{0}(u,\omega )=\mathrm{e}^{iu}X_{0}(\omega )$ for every $%
(u,\omega )\in \widetilde{\Omega }$, and for any $n\in {\mathbb{Z}}$, $%
\widetilde{X}_{n}=\widetilde{X}_{0}\circ {\widetilde{\theta }}^{n}$. Notice
that $(\widetilde{X}_{n})_{n\in {\mathbb{Z}}}$ is a stationary sequence of
complex random variables adapted to the non-decreasing filtration $(%
\widetilde{{\mathcal{F}}}_{n})$. Moreover $\mathrm{e}^{iu}\mathrm{e}%
^{int}X_{n}(\omega )=\widetilde{X}_{n}(u,\omega )$.

\medskip

\noindent \textbf{Proof of Proposition \ref{compactlilmart}}. Let $t\in
\lbrack 0,2\pi )$ be fixed. Let $\widetilde{D}(u)=\mathrm{e}^{iu}D$ and $%
\widetilde{D}_{k}=\widetilde{D}\circ \widetilde{\theta }^{k}$.

Let $\widetilde{d}_{k}=\big ({\mathcal{R}e}(\widetilde{D}_{k}),{\mathcal{I}m}%
(\widetilde{D}_{k})\big )^{\prime }$. Then 
\begin{equation*}
\widetilde{d}_{k}(u)=\left( 
\begin{array}{cc}
\cos u & -\sin u \\ 
\sin u & \cos u%
\end{array}%
\right) d_{k}(t)\,.
\end{equation*}%
Since the unit ball is invariant under rotations, the result will follow if we prove that for $%
\lambda $-a.e. $u\in \lbrack 0,2\pi ]$, the sequence $\big \{\sum_{k=1}^{n}%
\widetilde{d}_{k}(u)/\sqrt{2n\log \log n},n\geq 3\big \}$ has $%
\mbox{${\mathbb
P}$-a.s.}$ the ball $\{y\in {\mathbb{R}}^{2}\,:\,y^{\prime }\cdot y\leq
\Vert D_{0}\Vert _{2}^{2}/2\}$ as its set of limit points, or equivalently
(by Fubini's Theorem), if the sequence $\big \{\sum_{k=1}^{n}\widetilde{d}%
_{k}/\sqrt{2n\log \log n},n\geq 3\big \}$ has $%
\mbox{$\widetilde {\mathbb
P}$-a.s.}$ the ball $\{y\in {\mathbb{R}}^{2}\,:\,y^{\prime }\cdot y\leq
\Vert D\Vert _{2}^{2}/2\}$ as its set of limit points.

According to the almost sure analogue 
of the Cram\'{e}r-Wold device (see Sections 5.1 and 5.2 in Philipp \cite{Ph}%
), this will happen if we can prove that for any $x\in {\mathbb{R}}^{2}$
such that $x^{\prime }\cdot x=1$, 
\begin{equation}
\limsup_{n\rightarrow \infty }\frac{\sum_{k=1}^{n}x^{\prime }\cdot 
\widetilde{d}_{k}}{\sqrt{2n\log \log n}}=\frac{\Vert D\Vert _{2}^{2}}{2}%
\quad \widetilde{\mathbb{P}}-\text{a.s.}  \label{CWlil}
\end{equation}%
To prove it we shall apply Corollary 2 in Heyde and Scott \cite{HS} to the
stationary martingale differences $x^{\prime }\cdot \widetilde{d}\circ 
\widetilde{\theta }_{t}^{k}$. We have to verify 
\begin{equation}
\frac{1}{n}\sum_{k=1}^{n}(x^{\prime }\cdot \widetilde{d}_{k})^{2}\rightarrow 
\frac{\Vert D\Vert _{2}^{2}}{2}\qquad \widetilde{\mathbb{P}}-\text{a.s.}
\label{c0}
\end{equation}%
In order to understand this convergence it is convenient to write 
\begin{equation*}
D=A+iB\,.
\end{equation*}%
Therefore if $x=(a,b)^{\prime }$, 
\begin{gather*}
(x^{\prime }\cdot \widetilde{d}_{k})^{2}=\mathbb{(}a{\mathcal{R}e}(%
\widetilde{D}_{k})+b{\mathcal{I}m}(\widetilde{D}_{k}))^{2} \\
=(a\cos (u+kt)+b\sin (u+kt))^{2}(A^{2}\circ \theta ^{k})+(b\cos (u+kt)-a\sin
(u+kt))^{2}(B^{2}\circ \theta ^{k}) \\
+2(a\cos (u+kt)+b\sin (u+kt))(b\cos (u+kt)-a\sin (u+kt))(A\circ \theta
^{k})(B\circ \theta ^{k})\,.
\end{gather*}%
By using basic trigonometric formulas, it follows that if $x=(a,b)^{\prime }$
is such that $a^{2}+b^{2}=1$, 
\begin{gather*}
(x^{\prime }\cdot \widetilde{d}_{k})^{2}=\frac{(A^{2}+B^{2})\circ \theta ^{k}%
}{2}+\frac{(a^{2}-b^{2})\cos (2u+2kt)}{2}(A^{2}-B^{2})\circ \theta ^{k} \\
+\,ab\sin (2u+2kt)(A^{2}-B^{2})\circ \theta ^{k}+\,ab(\cos (2u+2kt)+\sin
(2u+2kt))(A\circ \theta ^{k})(B\circ \theta ^{k}) \\
+(b^{2}-a^{2})\sin (2u+2kt))(A\circ \theta ^{k})(B\circ \theta ^{k})\,.
\end{gather*}%
By Lemma \ref{LLNFourier} applied with $t_{0}=2t$, we derive that, for any $%
t\in (0,2\pi )\backslash \{\pi \}$ such that $\mathrm{e}^{-2it}$ is not an
eigenvalue of $\theta $ then, for all $u,$ 
\begin{equation}
\lim_{n\rightarrow \infty }\frac{1}{n}\sum_{k=1}^{n}\big ((x^{\prime }\cdot 
\widetilde{d}_{k}(u))^{2}-\frac{a^{2}+b^{2}}{2}(A^{2}+B^{2})\circ \theta ^{k}%
\big )=0\qquad \widetilde{\mathbb{P}}-\text{a.s.}  \label{c1}
\end{equation}%
Also by the ergodic theorem for $\theta $, 
\begin{equation}
\lim_{n\rightarrow \infty }\frac{1}{n}\sum_{k=1}^{n}(A^{2}+B^{2})\circ
\theta ^{k}=\mathbb{E}(|A|^{2}+|B|^{2})=\Vert D\Vert _{2}^{2}\qquad 
\widetilde{\mathbb{P}}-\text{a.s.}  \label{c2}
\end{equation}%
Gathering \eqref{c1} and \eqref{c2}, we get \eqref{c0}. This ends the proof
of Proposition \ref{compactlilmart}. \hfill $\square $

\medskip

\noindent \textbf{Proof of Proposition \ref{lmaasmartapprox}.}

\medskip

Let $\widetilde{D}_{0}(u,\cdot )=\mathrm{e}^{iu}\sum_{k\geq 0}P_{0}(\mathrm{e%
}^{ikt}X_{k})=\sum_{k\geq 0}\widetilde{P}_{0}(\widetilde{X}_{k})$ which is
defined for $\lambda $-a.e. $t\in \lbrack 0,2\pi )$. Write $\widetilde{S}%
_{n}=\sum_{k=1}^{n}\widetilde{X}_{k}$, $\widetilde{M}_{n}=\sum_{k=1}^{n}%
\widetilde{D}_{0}\circ \widetilde{\theta}^{k}$, and $\widetilde{R}_{n}=%
\widetilde{S}_{n}-\widetilde{M}_{n} $.

\medskip Denote by $\widetilde {\mathbb{E}} $ the expectation with respect
to $\widetilde {\mathbb{P}}$.

\medskip

The next lemma follows from Corollary 4.2 in Cuny \cite{Cuny}. Notice that
in \cite{Cuny}, complex-valued variables are allowed.

\begin{lemma}
\label{cuny} Assume that 
\begin{equation}  \label{condLIL}
\sum_{n\ge 1} \log n \frac{\widetilde {\mathbb{E}} (|\widetilde R_n|^2)}{n^2}%
<\infty\, .
\end{equation}
Then 
\begin{equation}  \label{asLIL}
\widetilde R_n=o(\sqrt{n\log \log n})\qquad 
\mbox{$\widetilde {\mathbb
P}$-a.s.}
\end{equation}
\end{lemma}

To prove that \eqref{approxmartlil} holds, it suffices to prove that for $%
\lambda $-a.e. $t\in \lbrack 0,2\pi )$, \eqref{asLIL} holds. According to
Lemma \ref{cuny} it suffices then to prove that \eqref{condLIL} is satisfied
for $\lambda $-a.e. $t\in \lbrack 0,2\pi )$. To this end, we first prove
that 
\begin{equation}
\int_{0}^{2\pi }\widetilde{{\mathbb{E}}}(|\widetilde{R}_{n}|^{2})dt
=2\sum_{k=1}^{n}\Vert {\mathbb{E}}_{0}(X_{k})\Vert _{2}^{2}\,.
\label{partest}
\end{equation}%
Indeed, for almost all $t \in [0,2\pi)$, 
\begin{equation}  \label{Rn}
\widetilde{R}_{n}(u,\omega )=\mathrm{e}^{iu}\sum_{k=1}^{n}\mathrm{e}^{ikt}{%
\mathbb{E}}_{0}(X_{k})(\omega )+\mathrm{e}^{iu}\sum_{k\geq n+1}\mathrm{e}%
^{ikt}({\mathbb{E}}_{0}(X_{k})(\omega )-{\mathbb{E}}_{n}(X_{k})(\omega ))\, .
\end{equation}
Whenever the R.H.S. below converges, the following identity holds: 
\begin{equation*}
\int_{0}^{2\pi }|\widetilde{R}_{n}(u,\omega )|^{2}dt=\sum_{k=1}^{n}({\mathbb{%
E}}_{0}(X_{k}))^{2}(\omega )+\sum_{k\geq n+1}({\mathbb{E}}_{0}(X_{k})(\omega
)-{\mathbb{E}}_{n}(X_{k})(\omega ))^{2}\,.
\end{equation*}%
Then, using that ${\mathbb{E}}(({\mathbb{E}}_{0}(X_{k})-{\mathbb{E}}%
_{n}(X_{k}))^{2})=\Vert {\mathbb{E}}_{0}(X_{k-n})\Vert _{2}^{2}-\Vert {%
\mathbb{E}}_{0}(X_{k})\Vert _{2}^{2}$, and the fact that under %
\eqref{condFM2}, \eqref{regular} holds, we obtain \eqref{partest}. Using %
\eqref{partest}, we see that under \eqref{condFM2}, for $\lambda $-a.e. $%
t\in \lbrack 0,2\pi )$, condition \eqref{condLIL} holds. This ends the proof
of \eqref{approxmartlil} and then of the proposition. \hfill $\square $

\bigskip

\noindent \textbf{Proof of Theorem \ref{thmlilmoment2delta}.} According to
the proof of Theorem \ref{thmlilmoment2}, it suffices to prove that under
the conditions of Theorem \ref{thmlilmoment2delta}, the almost sure
convergence \eqref{approxmartlil} holds for almost all $t\in \lbrack 0,2\pi
) $. With this aim, we shall use truncation arguments. Given $\gamma >0$ and 
$r\geq 0$, we set for any integer $\ell $, 
\begin{equation*}
\overline{X}_{\ell ,r}:=X_{\ell }\mathbf{1}_{\{|X_{\ell }|\leq 2^{\gamma
r}\}}-{\mathbb{E}}(X_{\ell }\mathbf{1}_{\{|X_{\ell }|\leq 2^{\gamma r}\}})
\end{equation*}%
and 
\begin{equation*}
\overline{D}_{\ell ,r}(t):=\mathrm{e}^{i\ell t}\sum_{k\geq 0}\mathrm{e}^{ikt}%
\big (P_{0}(\overline{X}_{k,r})\big )\circ \theta ^{\ell }\,.
\end{equation*}%
We know that for almost all $t\in [0,2\pi )$, $\overline{D}_{\ell ,r}(t)$ is
defined $\mbox{${\mathbb P}$-a.s.}$ and in ${\mathbb{L}}^{2}({\mathbb{P}})$.

We define non stationary sequences $(\overline{X}_{\ell })_{\ell \geq 1}$
and $(\overline{D}_{\ell }(t))_{\ell \geq 1}$ as follows: for every $r\in {%
\mathbb{N}}$ and every $\ell \in \{2^{r},...,2^{r+1}-1\}$, 
\begin{equation}
\overline{X}_{\ell }:=\overline{X}_{\ell ,r}\,,\ \overline{D}_{\ell }(t):=%
\overline{D}_{\ell ,r}(t)\,.  \label{defbar}
\end{equation}%
Let also 
\begin{equation}
X_{\ell }^{\ast }=X_{\ell }-\overline{X}_{\ell }\ \text{ and }\ D_{\ell
}^{\ast }(t)=D_{\ell }(t)-\overline{D}_{\ell }(t)\,.  \label{deftilde}
\end{equation}

\begin{lemma}
\label{lmatruncation1} Assume that ${\mathbb{E}}\Big (\frac{X_{0}^{2}L(X_0)}{%
L(L(X_{0}))}\Big )<\infty $. Then, for a.e. $t\in [0,2\pi )$,%
\begin{equation}
\sum_{n\geq 3}\frac{X_{n}^{\ast }}{\sqrt{n\log \log n}}\mathrm{e}^{int}\quad 
\text{converges }\mathbb{P}\text{-a.s.}  \label{xstar}
\end{equation}%
In particular, by Kronecker's lemma, for a.e. $t\in \lbrack 0,2\pi )$, $%
\frac{\sum_{k=1}^{n}\mathrm{e}^{ikt}X_{k}^{\ast }}{\sqrt{n\log \log n}}%
\rightarrow 0$ $\mbox{${\mathbb P}$-a.s.}$
\end{lemma}

\noindent \textbf{Proof.} By Carleson's theorem \cite{Ca}, in order to
establish (\ref{xstar}) it suffices to prove that%
\begin{equation*}
\sum_{n\geq 3}\frac{(X_{n}^{\ast })^{2}}{n\log \log n}<\infty \qquad 
\mbox{${\mathbb
P}$-a.s.}
\end{equation*}%
This is true because 
\begin{align}
\sum_{n\geq 4}\frac{{\mathbb{E}}((X_{n}^{\ast })^{2})}{n\log \log n}&
=\sum_{r\geq 2}\sum_{\ell =2^{r}}^{2^{r+1}-1}\frac{{\mathbb{E}}((X_{\ell }-%
\overline{X}_{\ell ,r})^{2})}{\ell \log \log \ell }\leq 4\sum_{r\geq 2}{%
\mathbb{E}}(X_{0}^{2}\mathbf{1}_{\{|X_{0}|>2^{\gamma r}\}})\sum_{\ell
=2^{r}}^{2^{r+1}-1}\frac{1}{\ell \log \log \ell }  \notag
\label{p1lma1trunc} \\
& \ll \sum_{r\geq 2}\frac{1}{\log r}{\mathbb{E}}(X_{0}^{2}\mathbf{1}%
_{\{|X_{0}|>2^{\gamma r}\}})\ll {\mathbb{E}}\Big (\frac{X_{0}^{2}L(X_{0})}{%
L(L(X_{0}))}\Big )<\infty \,,
\end{align}%
where we used Fubini in the last step and the notation $a\ll b$
means there is a universal constant $C>0$ such that $a<Cb$. \hfill $\square $

\begin{lemma}
\label{lmatruncation2} Assume that ${\mathbb{E}}\Big (\frac{X_{0}^{2}L(X_{0})%
}{L(L(X_{0}))}\Big )<\infty $. Then, for a.e. $t\in \lbrack 0,2\pi )$, 
\begin{equation*}
\frac{\sum_{k=1}^{n}D_{k}^{\ast }(t)}{\sqrt{n\log \log n}}\rightarrow 0\quad %
\mbox{${\mathbb P}$-a.s.}
\end{equation*}
\end{lemma}

\noindent \textbf{Proof.} For almost all $t\in \lbrack 0,2\pi )$, $(D_{\ell
}^{\ast }(t))_{\ell \geq 1}$ is a sequence of martingale differences in ${%
\mathbb{L}}^{2}({\mathbb{P}})$. Hence using the Doob-Kolmogorov maximal
inequality, we infer that the lemma will be established provided that 
\begin{equation}
\sum_{k\geq 3}\int_{0}^{2\pi }\frac{\Vert D_{k}^{\ast }(t)\Vert _{2}^{2}}{%
k\log \log k}dt<\infty \,.  \label{p1lma2trun}
\end{equation}%
To prove it we use simple algebra and the projection's orthogonality, as
follows:%
\begin{gather*}
\sum_{k\geq 4}\int_{0}^{2\pi }\frac{\Vert D_{k}^{\ast }(t)\Vert _{2}^{2}}{%
k\log \log k}dt=\sum_{r\geq 2}\sum_{\ell =2^{r}}^{2^{r+1}-1}\int_{0}^{2\pi }%
\frac{\Vert D_{\ell }(t)-\overline{D}_{\ell ,r}(t)\Vert _{2}^{2}}{\ell \log
\log \ell }dt \\
\leq \sum_{r\geq 2}\frac{1}{2^{r}\log r}\sum_{\ell
=2^{r}}^{2^{r+1}-1}\int_{0}^{2\pi }\frac{\Vert D_{\ell }(t)-\overline{D}%
_{\ell }(t)\Vert _{2}^{2}}{\ell \log \log \ell }dt \\
\leq 2\pi \sum_{r\geq 2}\frac{1}{2^{r}\log r}\sum_{\ell
=2^{r}}^{2^{r+1}-1}\sum_{k\geq 0}\Vert \big (P_{0}(X_{k}-\overline{X}_{k,r})%
\big )\circ \theta ^{\ell }\Vert _{2}^{2} \\
=2\pi \sum_{r\geq 2}\frac{1}{2^{r}\log r}\sum_{\ell
=2^{r}}^{2^{r+1}-1}\sum_{k\geq 0}\Vert \big (P_{-k}(X_{0}\mathbf{1}%
_{\{|X_{0}|>2^{\gamma r}\}})\big )\circ \theta ^{\ell +k}\Vert _{2}^{2} \\
=2\pi \sum_{r\geq 2}\frac{1}{\log r}{\mathbb{E}}\big (X_{0}\mathbf{1}%
_{\{|X_{0}|>2^{\gamma r}\}}-{\mathbb{E}}(X_{0}\mathbf{1}_{\{|X_{0}|>2^{%
\gamma r}\}}|{\mathcal{F}}_{-\infty })\big )^{2}\leq 2\pi \sum_{r\geq 2}%
\frac{1}{\log r}\Vert X_{0}\mathbf{1}_{\{|X_{0}|>2^{\gamma r}\}}\Vert
_{2}^{2}\,.
\end{gather*}%
Next, using Fubini's theorem as done in \eqref{p1lma1trunc}, %
\eqref{p1lma2trun} follows. \hfill $\square $

\medskip

From Lemmas \ref{lmatruncation1} and \ref{lmatruncation2}, we then deduce
that if ${\mathbb{E}}\Big (\frac{X_{0}^{2}L(X_{0})}{L(L(X_{0}))}\Big )%
<\infty $, then, for a.e. $t\in \lbrack 0,2\pi )$, 
\begin{equation*}
\frac{\sum_{k=1}^{n}\big (\mathrm{e}^{ikt}X_{k}^{\ast }-D_{k}^{\ast }(t)\big
)}{\sqrt{n\log \log n}}\rightarrow 0\quad \mbox{${\mathbb P}$-a.s.}
\end{equation*}%
Therefore, to prove that the almost sure convergence \eqref{approxmartlil}
holds for almost all $t\in \lbrack 0,2\pi )$ (and then the theorem) it
suffices to prove that for almost all $t\in \lbrack 0,2\pi )$, 
\begin{equation}
\frac{|\overline{S}_{n}(t)-\overline{M}_{n}(t)|}{\sqrt{n\log \log n}}%
\rightarrow 0\ \ \mbox{${\mathbb P}$-a.s.}  \label{approxmartlil2}
\end{equation}%
where $\overline{S}_{n}(t)=\sum_{j=1}^{n}\mathrm{e}^{ijt}\overline{X}_{j}$
and $\overline{M}_{n}(t)=\sum_{j=1}^{n}\overline{D}_{j}(t)$ where the $%
\overline{X}_{j}$'s and $\overline{D}_{j}(t)$'s are defined in \eqref{defbar}%
. Let 
\begin{equation*}
\overline{R}_{n}(t)=\overline{S}_{n}(t)-\overline{M}_{n}(t)\,,
\end{equation*}%
and for any $r\in {\mathbb{N}}$, let 
\begin{equation*}
A_{r}(t):=\sup_{0\leq k\leq 2^{r}-1}|\overline{R}_{k+2^{r}}(t)-\overline{R}%
_{2^{r}-1}(t)|\,.
\end{equation*}%
Let $N\in {\mathbb{N}}^{\ast }$ and let $k\in ]1,2^{N}]$. We first notice
that $A_{r}(t)\geq |\sum_{\ell =2^{r}}^{2^{r+1}-1}\mathrm{e}^{i\ell t}(%
\overline{X}_{\ell }-\overline{D}_{\ell }(t))|$, so if $K$ is the integer
such that $2^{K-1}\leq k\leq 2^{K}-1$, then 
\begin{equation*}
\big |\overline{R}_{k}(t)\big |\leq \sum_{r=0}^{K-1}A_{r}(t)\,.
\end{equation*}%
Consequently since $K\leq N$, 
\begin{equation*}
\sup_{1\leq k\leq 2^{N}}\big |\overline{R}_{k}(t)\big |\leq
\sum_{r=0}^{N-1}A_{r}(t)\,.
\end{equation*}%
Therefore, \eqref{approxmartlil2} will follow if we can prove that for
almost all $t\in \lbrack 0,2\pi )$, 
\begin{equation*}
\sup_{0\leq k\leq 2^{r}-1}|\overline{R}_{k+2^{r}}(t)-\overline{R}%
_{2^{r}-1}(t)|=o\big (2^{r/2}\cdot (\log r)^{1/2}\big )\quad 
\mbox{${\mathbb
P}$-a.s.}\,,
\end{equation*}%
which will be true if we can prove that 
\begin{equation}
\sum_{r\geq 0}\frac{1}{2^{r}\log r}\int_{0}^{2\pi }{\mathbb{E}}%
[\max_{2^{r}\leq k\leq 2^{r+1}-1}|\overline{R}_{k}(t)-\overline{R}%
_{2^{r}-1}(t)|^{2}]dt<\infty \,.  \label{conv1thm2}
\end{equation}%
Notice that for any integer $k$ in $[2^{r},2^{r+1}-1]$, 
\begin{equation*}
\overline{R}_{k}(t)-\overline{R}_{2^{r}-1}(t)=\mathrm{e}^{i(2^{r}-1)t}\Big (%
\sum_{\ell =1}^{k-2^{r}+1}\mathrm{e}^{i\ell t}\big (\overline{X}_{0,r}-%
\overline{D}_{0,r}(t)\big )\circ \theta ^{\ell }\Big )\circ \theta
^{2^{r}-1}\,.
\end{equation*}%
Therefore, by stationarity proving \eqref{conv1thm2} amounts to prove that 
\begin{equation}
\sum_{r\geq 0}\frac{1}{2^{r}\log r}\int_{0}^{2\pi }{\mathbb{E}}\big (%
\max_{1\leq k\leq 2^{r}}\big |\sum_{\ell =1}^{k}(\overline{X}_{\ell ,r}-%
\overline{D}_{\ell ,r}(t))\big |^{2}\big )dt<\infty \,,  \label{conv1thm2bis}
\end{equation}%
where $\overline{X}_{\ell ,r}(t)=\mathrm{e}^{i\ell t}\overline{X}_{\ell ,r}$%
. Using Lemma \ref{max2} given in the Appendix with $M=2^{\gamma r}$, we get that for any integer $s>1$, 
\begin{align*}
(2\pi )^{-1}\int_{0}^{2\pi }& {\mathbb{E}}\big (\max_{1\leq k\leq 2^{r}}\big
|\sum_{\ell =1}^{k}(\overline{X}_{\ell ,r}-\overline{D}_{\ell ,r}(t))\big |%
^{2}\big )dt \\
& \leq 24\times 2^{r}\Vert \mathbb{E}_{-s}(X_{0})\Vert ^{2}+24\times
2^{r}\Vert X_{0}\mathbf{1}_{|X_{0}|>2^{\gamma r}}\Vert
^{2}+12s^{2}2^{2\gamma r}\,.
\end{align*}%
To prove \eqref{conv1thm2bis} and then to end the proof of the theorem, we
select $\gamma <1/4$ and use the above inequality with $s=[2^{\gamma r}]+1$.
Using Fubini's theorem as done in \eqref{p1lma1trunc}, we infer that %
\eqref{conv1thm2bis} will be established provided that 
\begin{equation}
\sum_{r\geq 2}\frac{1}{\log r}\Vert {\mathbb{E}}_{-[2^{\gamma
r}]}(X_{0})\Vert _{2}^{2}<\infty \,.  \label{convseries2}
\end{equation}%
This convergence follows from condition \eqref{condFM} by using the fact $%
(\Vert {\mathbb{E}}_{-n}(X_{0})\Vert _{2}^{2})_{n\geq 1}$ is decreasing and
by noticing that by the usual comparison between the series and the
integrals, for any non-increasing and positive function $h$ on ${\mathbb{R}}%
^{+}$ and any positive $\gamma $, 
\begin{equation*}
\sum_{n\geq 1}n^{-1}h(n^{\gamma })<\infty \text{ if and only if }\sum_{n\geq
1}n^{-1}h(n)<\infty \,.
\end{equation*}%
and that \eqref{convseries2} is equivalent to $\sum_{n\geq 3}\frac{1}{n(\log
\log n)}\Vert {\mathbb{E}}_{-[n^{\gamma }]}(X_{0})\Vert _{2}^{2}<\infty $.
\hfill $\square $

\bigskip

\noindent \textbf{Proof of Theorem \ref{thm3}.} We divide the proof of this
theorem in two parts.

\medskip

1. \textbf{Proof of the continuity of $f$ and of relation (\ref{defgt}).}
Let $(c_{n})_{n\in \mathbb{Z}}$ denote the Fourier coefficients of $f$, i.e. 
$c_{n}:=\mathbb{E}(X_{0}X_{n})$. Then, the Fourier coefficients of $(1-%
\mathrm{e}^{it})f(t)$ are $(c_{n}-c_{n+1})_{n\in \mathbb{Z}}$ and the
Fourier coefficients of $h(t):=|1-\mathrm{e}^{it}|^{2}f(t)$ are $%
(b_{n})_{n\in \mathbb{Z}}$ with $b_{n}=2c_{n}-c_{n+1}-c_{n-1},$ $n\in 
\mathbb{Z}$.

One can easily see that $h$ is the spectral density associated with the
stationary process $(Z_{n})_{n\in \mathbb{Z}}:=(X_{n}-X_{n-1})_{n\in \mathbb{%
Z}}$, i.e. $b_{n}=\mathbb{E}(Z_{0}Z_{n})$. Hence for $n\geq 0$, 
\begin{equation*}
|b_{n}|=|\mathbb{E}(Z_{0}Z_{n})|=|\sum_{k\geq 0}\mathbb{E}%
(P_{-k}(Z_{0})P_{-k}(Z_{n}))|\leq \sum_{k\geq 0}\Vert P_{-k}(Z_{0})\Vert
_{2}\Vert P_{-k-n}(Z_{0})\Vert _{2}\,.
\end{equation*}%
Therefore, 
\begin{equation*}
\sum_{n\in \mathbb{Z}\backslash \{0\}}|b_{n}|=2\sum_{n\geq 1}|b_{n}|\leq
2\sum_{n\geq 0}\big(\Vert P_{-n}(Z_{0})\Vert _{2}\big)^{2}\,.
\end{equation*}%
By (\ref{condproj}) it follows that $(b_{n})_{n\in \mathbb{Z}}$ is
absolutely summable. Therefore, by well known results on spectral density,
(see for instance Bradley \cite{Br}, Ch 8 and 9) $h$ must be continuous and bounded on $[0,2\pi ]$, which in turn
implies that $f$ is continuous on $(0,2\pi )$.

We prove now that (\ref{defgt}) holds for every $t\in (0,2\pi )$. With this
aim, it suffices to show that, for every $t\in (0,2\pi )$, $|1-\mathrm{e}%
^{it}|^{2}\mathbb{E}(|S_{n}(t)|^{2})/n\rightarrow h(t)$.

Define $T_{n}(t):=\sum_{k=1}^{n}\mathrm{e}^{ikt}Z_{k}$. It is easy to see,
using the fact that $c_{n}\rightarrow 0$ as $n\rightarrow \pm \infty $, that 
\begin{equation*}
\frac{\mathbb{E}(|T_{n}(t)|^{2})}{n}=\frac{|1-\mathrm{e}^{it}|^{2}\mathbb{E}%
(|S_{n}(t)|^{2})}{n}+o(1)\,,
\end{equation*}%
(the little $o$ is uniform in $t\in \lbrack 0,2\pi ])$. Now, by (4.6) of 
\cite{PeWu}, $(\mathbb{E}(|T_{n}(t)|^{2}/n)_{n\geq 0}$ is nothing else but
the Ces\`{a}ro averages of the partial sums of the Fourier series associated
with $h $, hence it converges to $h(t)$ by Fejer's theorem.

\medskip

2. \textbf{End of the proof.}

By assumption (\ref{condproj}), we have 
\begin{equation}
\sum_{k\geq 0}|P_{0}(X_{k}-X_{k+1})|\qquad \mbox{converges  in $\LL^2$}.
\label{assump}
\end{equation}%
Let $t\in (0,2\pi )$ be \emph{fixed}. Using that $P_{0}(X_{-1})=0$, we
obtain 
\begin{gather*}
\sum_{m=0}^{k}\mathrm{e}^{imt}P_{0}(X_{m}-X_{m-1})=\sum_{m=0}^{k}\mathrm{e}%
^{imt}P_{0}(X_{m})-\sum_{m=0}^{k-1}\mathrm{e}^{i(m+1)}P_{0}(X_{m}) \\
=(1-\mathrm{e}^{it})\sum_{m=0}^{k}\mathrm{e}^{imt}P_{0}(X_{m})+\mathrm{e}%
^{i(k+1)t}P_{0}(X_{k})\,.
\end{gather*}%
Since $\Vert P_{0}(X_{k})\Vert _{2}\rightarrow 0$, by (\ref{assump}), we see %
\label{some little changes} that the series $\sum_{m=0}^{k}\mathrm{e}%
^{imt}P_{0}(X_{m})$ converges in ${\mathbb{L}}^2$ as $k\rightarrow \infty $.
Hence defining $D_0(t)$ by \eqref{defD0alpha}, it follows that $(D_{0} (t)
\circ \theta^{\ell} , \ell \in {\mathbb{Z}} )$ is a stationary sequence of
martingale differences in ${\mathbb{L}}^2$ adapted to $( {\mathcal{F}}%
_{\ell} , \ell \in {\mathbb{Z}} )$. Hence the theorem will follow by
Proposition \ref{compactlilmart}, if we can prove that $|S_{n}(t)-M_{n}(t)|/%
\sqrt{n\log \log n}\rightarrow 0$ $\mathbb{P}$-a.s. where $M_n(t) =
\sum_{k=1}^n \mathrm{e}^{ikt} D_0(t) \circ \theta^k$. With this aim, we
first notice that 
\begin{equation*}
(1-\mathrm{e}^{it})D_{0}(t)=F_{0}(t)\qquad \mbox{where}\qquad
F_{0}(t)=\sum_{m\geq 0}\mathrm{e}^{imt}P_{0}(X_{m}-X_{m-1})\,.
\end{equation*}%
Hence, writing $F_{k}(t)=F_{0}(t)\circ \theta ^{k}$, we obtain the
representation 
\begin{equation*}
(1-\mathrm{e}^{it})(S_{n}(t)-M_{n}(t))=\sum_{k=0}^{n-1}\mathrm{e}%
^{ikt}(Z_{k}-F_{k}(t))\,,
\end{equation*}%
where $Z_{k}=X_{k}-X_{k-1}$. Therefore, the proof of the theorem will be
complete if we can show that 
\begin{equation}  \label{modifconvps}
\big|\sum_{k=0}^{n-1}\mathrm{e}^{ikt}(Z_{k}-F_{k}(t))\big|/\sqrt{n\log \log n%
}\rightarrow 0\qquad \mbox{${\mathbb P}$-a.s.}
\end{equation}
To prove this almost sure convergence, we shall work on the product space $(%
\widetilde{\Omega },\widetilde{{\mathcal{F}}},\widetilde{\mathbb{P}})$
introduced in the proof of Theorem \ref{thmlilmoment2}. Recall that $%
\widetilde{\theta}$ has been defined in \eqref{thetatilde}, $\widetilde {%
\mathcal{F}}_n:= {\mathcal{B} ([0,2\pi])}\otimes {\mathcal{F}}_n$ and $%
\widetilde {\mathbb{E}}$ stands for the expectation under $\widetilde{%
\mathbb{P}}$.

Define $\widetilde{Z}_{0}(u,\omega ):=\mathrm{e}^{iu}Z_{0}(\omega )$ and $%
\widetilde{Z}_{k}:=\widetilde{Z}_{0}\circ \widetilde{\theta}^{k}$.
Similarly, define $\widetilde{F}_0(u,\omega )=\mathrm{e}^{iu}F_{0}(t)(\omega
)$ and $\widetilde{F}_{k}=\widetilde{F}_0 \circ \widetilde{\theta}^{k}$. Let 
$\widetilde{P}_{0}( \cdot ) = \widetilde{{\mathbb{E}}}( \cdot | \widetilde {%
\mathcal{F}}_0) - \widetilde{{\mathbb{E}}}( \cdot | \widetilde {\mathcal{F}}%
_{-1})$. Note that $\widetilde{F}_{0}=\sum_{k\geq 0}\widetilde{P}_{0}(%
\widetilde{Z}_{k})=\mathrm{e}^{iu}\sum_{k\geq 0}\mathrm{e}%
^{ikt}P_{0}(X_{k}-X_{k-1})$.

By assumption (\ref{condproj}), we have 
\begin{equation*}
\sum_{n\geq 0}\Vert \widetilde{P}_{0}(\widetilde{Z}_{n})\Vert _{2, 
\widetilde{\mathbb{P}}}<\infty \, ,
\end{equation*}
where $\Vert \cdot \Vert _{2, \widetilde{\mathbb{P}}}$ is the ${\mathbb{L}}%
^2 $ norm with respect to $\widetilde{\mathbb{P}}$.

Therefore, by Theorem 2.7 of Cuny \cite{CunyLIL} (see (21) of \cite{CunyLIL}%
), identifying $\mathbb{C}$ with $\mathbb{R}^{2}$, we obtain that 
\begin{equation*}
\big|\sum_{k=0}^{n-1}(\widetilde{Z}_{k}-\widetilde{F}_{k})\big|/\sqrt{n\log
\log n}\rightarrow 0\qquad \mbox{$\widetilde \P$-a.s.}
\end{equation*}%
Now, $\widetilde{Z}_{k}(u,.)=\mathrm{e}^{iu}\mathrm{e}^{ikt}Z_{k}$ and $%
\widetilde{F}_{k}(u,.)=\mathrm{e}^{iu}\mathrm{e}^{ikt}F_{k}$, hence %
\eqref{modifconvps} follows. \hfill $\square $

\bigskip

\noindent \textbf{Proof of Theorem \ref{thm4}.} Define an operator $R_{t}$
on $\mathbb{L}^{2}(\Omega ,{\mathcal{F}}_{0},\mathbb{P})$ by $R_{t}(Y):=%
\mathrm{e}^{it}{\mathbb{E}}_{0}(Y\circ \theta )\,.$ Note that for every $%
n\geq 0$, $R_{t}^{n}(Y)=\mathrm{e}^{int}{\mathbb{E}}_{0}(Y\circ \theta ^{n})$%
. Hence by assumption $\sup_{n\geq 1}\Vert
\sum_{k=0}^{n}R_{t}^{k}(X_{0})\Vert _{2}<\infty \,.$ By Browder \cite[Lemma 5%
]{Browder}, there exists $Z_{0}=Z_{0}(t)\in {\mathbb{L}}^{2}(\Omega ,{%
\mathcal{F}}_{0},\mathbb{P})$ such that 
\begin{equation}
X_{0}=Z_{0}-R_{t}(Z_{0})=Z_{0}-\mathrm{e}^{it}{\mathbb{E}}_{0}(Z_{1})\,.
\label{Poisson}
\end{equation}%
Now we denote $Z_{k}=Z_{0}\circ \theta ^{k}$. Note that $(Z_{k})_{k}$ is a
stationary sequence, $R_{t}(Z_{0})=\mathrm{e}^{it}{\mathbb{E}}_{0}(Z_{1})\ $%
and we have the decomposition: 
\begin{equation*}
X_{0}=Z_{0}-{\mathbb{E}}_{-1}(Z_{0})+{\mathbb{E}}_{-1}(Z_{0})-\mathrm{e}^{it}%
{\mathbb{E}}_{0}(Z_{1})\,.
\end{equation*}%
Denote the martingale difference $D_{0}(t)=Z_{0}-{\mathbb{E}}%
_{-1}(Z_{0})=P_{0}(Z_{0})$. So,%
\begin{eqnarray*}
S_{n}(t) &=&\sum_{k=1}^{n}\mathrm{e}^{itk}D_{0}(t)\circ \theta
^{k}+\sum_{k=1}^{n}(\mathrm{e}^{itk}{\mathbb{E}}_{k-1}(Z_{k})-\mathrm{e}%
^{it(k+1)}{\mathbb{E}}_{k}(Z_{k+1})) \\
&=&\,\sum_{k=1}^{n}\mathrm{e}^{itk}D_{0}(t)\circ \theta ^{k}+\mathrm{e}^{it}{%
\mathbb{E}}_{0}(Z_{1})-\mathrm{e}^{it(n+1)}{\mathbb{E}}_{n}(Z_{n+1})\,.
\end{eqnarray*}%
By the Borel-Cantelli lemma,%
\begin{equation*}
|{\mathbb{E}}_{n}(Z_{n+1})|/\sqrt{n}=|{\mathbb{E}}_{-1}(Z_{t})|\circ \theta
^{n+1}/\sqrt{n}\rightarrow 0\qquad \mbox{${\mathbb P}$-a.s.}
\end{equation*}%
Therefore we have the following martingale approximation: 
\begin{equation*}
\frac{1}{n}\big (S_{n}(t)-\sum_{k=1}^{n}\mathrm{e}^{itk}D_{0}(t)\circ \theta
^{k}\big )\rightarrow 0\text{ a.s. and in }\mathbb{L}^{2}\text{ .}
\end{equation*}%
Hence, since $\mathrm{e}^{-2it}$ is not an eigenvalue of $\theta $, the
proposition follows from Proposition \ref{compactlilmart} with $%
E(|D_{0}|^{2})=E(|Z_{0}-{\mathbb{E}}_{-1}(Z_{0})|^{2})=\sigma _{t}^{2}.$

It is convenient to express $\sigma _{t}$ in terms of the original
variables. With this aim notice that by equation (\ref{Poisson}) we obtain%
\begin{gather*}
\sum_{k=0}^{n}\mathrm{e}^{itk}P_{0}(X_{k})=\sum_{k=0}^{n}\mathrm{e}%
^{itk}P_{0}(Z_{k})-\sum_{k=0}^{n}\mathrm{e}^{it(k+1)}P_{0}(Z_{k+1}) \\
=P_{0}(Z_{0})-\mathrm{e}^{it(n+1)}P_{0}(Z_{n+1})=D_{0}(t)-\mathrm{e}%
^{it(n+1)}P_{0}(Z_{n+1})\,.
\end{gather*}%
Since 
\begin{equation*}
\Vert P_{0}(Z_{n+1})\Vert _{2}^{2}=\Vert \mathbb{E}_{0}(Z_{n+1})\Vert
_{2}^{2}-\Vert \mathbb{E}_{-1}(Z_{n+1})\Vert _{2}^{2}=\Vert \mathbb{E}%
_{-n-1}(Z_{0})\Vert _{2}^{2}-\Vert \mathbb{E}_{-n-2}(Z_{0})\Vert
_{2}^{2}\rightarrow 0\,,
\end{equation*}%
we obtain 
\begin{equation*}
\sum_{k=0}^{n}\mathrm{e}^{itk}P_{0}(X_{k})\rightarrow D_{0}(t)\text{ in }{%
\mathbb{L}}^{2}\,.
\end{equation*}%
This shows that, for this case, the representation (\ref{defD0alpha}) holds
for all $t\in \lbrack 0,2\pi ]$ such that \eqref{poisson} is satisfied.
\hfill $\square $

\section{Examples}

\label{sectappl}

\subsection{Linear processes.}

\label{sectionLP}

Let us consider the following linear process $(X_{k})_{k\in {\mathbb{Z}}}$
defined by $X_{k}=\sum_{j\geq 0}a_{j}\varepsilon _{k-j}$ where $(\varepsilon
_{k})_{k\in {\mathbb{Z}}}$ is a sequence of iid real random variables in ${%
\mathbb{L}}^{2}$ and $(a_{k})_{k\in {\mathbb{Z}}}$ is a sequence of reals in 
$\ell ^{2}$. Taking ${\mathcal{F}}_{0}=\sigma (\varepsilon _{k},k\leq 0)$,
it follows that $P_{0}(X_{i})=a_{i}\varepsilon _{0}$. Therefore %
\eqref{condproj} is reduced to 
\begin{equation}
\sum_{n\geq 3}|a_{n}-a_{n+1}|<\infty \,.  \label{condlinear}
\end{equation}%
Hence, because ${\mathcal{F}}_{-\infty }$ is trivial, we conclude, by
Corollary \ref{cor5}, that the conclusions of Theorem \ref{thm3} hold for
all $t\in (0,\pi )\cup (\pi ,2\pi )$ as soon as \eqref{condlinear} is
satisfied. Let us mention that when $a_{n}$ is decreasing (\ref{condlinear})
is always satisfied.

\subsection{Functions of linear processes}

In this section, we shall focus on functions of real-valued linear
processes. Define 
\begin{equation}
X_{k}=h\Big(\sum_{i\geq 0}a_{i}\varepsilon _{k-i}\Big)-{\mathbb{E}}\Big(h%
\Big(\sum_{i\geq 0}a_{i}\varepsilon _{k-i}\Big)\Big)\, ,  \label{def2suite}
\end{equation}
where $(a_{i})_{i\in {{\mathbf{Z}}}}$ be a sequence of real numbers in $\ell
^{2}$ and $(\varepsilon _{i})_{i\in \mathbb{Z}}$ is a sequence of iid random
variables in ${\mathbb{L}}^{2}$. We shall give sufficient conditions in
terms of the regularity of the function $h$, for $I_n(t)$ to satisfy a law
of the iterated logarithm as described in Theorem \ref{thmlilmoment2}.

Denote by $w_{h}(.,M)$ the modulus of continuity of the function $h$ on the
interval $[-M,M]$, that is 
\begin{equation*}
w_{h}(u,M)=\sup \{|h(x)-h(y)|,|x-y|\leq u,|x|\leq M,|y|\leq M\}\,.
\end{equation*}

\begin{corollary}
\label{corlin} Assume that $h$ is $\gamma $-H\"{o}lder on any compact set,
with $w_{h}(u,M)\leq Cu^{\gamma }M^{\alpha }$, for some $C>0$, $\gamma \in
]0,1]$ and $\alpha \geq 0$. Assume that 
\begin{equation}
\sum_{k\geq 2}(\log k)^{2}|a_{k}|^{2\gamma }<\infty \text{ and }{\mathbb{E}}%
(|\varepsilon _{0}|^{2\vee (2\alpha +2\gamma )})<\infty \,.
\label{condcorflin}
\end{equation}%
Then the conclusions of Theorem \ref{thmlilmoment2} hold with $(X_{k})_{k\in 
{\mathbf{Z}}}$ defined by (\ref{def2suite}).
\end{corollary}

\noindent\textbf{Proof.} We shall apply Theorem \ref{thmlilmoment2} by
taking ${\mathcal{F}}_k = \sigma( \varepsilon_{\ell} , \ell \leq k)$. Since $%
X_0$ is regular, $\Vert {\mathbb{E}}_0 ( X_k ) \Vert_2^2 = \sum_{\ell \geq
k} \Vert P_{-\ell} (X_0) \Vert_2^2$. Therefore \eqref{condFM2} is equivalent
to 
\begin{equation}  \label{condFM2bis}
\sum_{\ell \geq 2} (\log \ell)^2 \Vert P_{0} (X_{\ell}) \Vert_2^2<\infty \, .
\end{equation}
Let $\varepsilon^{\prime }$ be an independent copy of $\varepsilon$, and
denote by ${\mathbb{E}}_{\varepsilon}(\cdot)$ the conditional expectation
with respect to $\varepsilon$. Clearly 
\begin{equation*}  \label{bof}
P_0(X_k)= {\mathbb{E}}_{\varepsilon}\Big(h\Big(\sum_{i=0 }^{k-1} a_i
\varepsilon^{\prime }_{k-i} +a_k \varepsilon_0 + \sum_{i >k} a_i
\varepsilon_{k-i}\Big) - h\Big(\sum_{i=0 }^{k-1} a_i \varepsilon^{\prime
}_{k-i} +a_k \varepsilon^{\prime }_0 + \sum_{i >k} a_i \varepsilon_{k-i}%
\Big) \Big) \, .
\end{equation*}
Since $w_{h}(u_1+u_2, M) \leq w_{h}(u_1, M) + w_{h}(u_2, M)$, it follows
that 
\begin{equation}  \label{majpo}
|P_0(X_k)| \leq {\mathbb{E}}_{\varepsilon}\Big( 2\|X_0\|_\infty \wedge \Big(%
w_{h}\left( |a_k||\varepsilon_0|, |Y_1| \vee |Y_2|\right) + w_{h}\left(
|a_k||\varepsilon^{\prime }_0|, |Y_1|\vee |Y_2|\right) \Big) \Big)\, ,
\end{equation}
where $Y_1= \sum_{i=0}^k a_i \varepsilon^{\prime }_{k-i} + \sum_{i >k} a_i
\varepsilon_{k-i}$ and $Y_2= \sum_{i =0}^{k-1} a_i \varepsilon^{\prime
}_{k-i} + \sum_{i \geq k} a_i \varepsilon_{k-i}$. Noting that $%
(\varepsilon_0, |Y_1|\vee|Y_2|)$ and $(\varepsilon_0^{\prime },
|Y_1|\vee|Y_2|)$ are both distributed as $(\varepsilon_0, M_k)$, where $M_k=
\max \Big \{ \Big |\sum_{i \geq 0} a_i \varepsilon^{\prime }_i\Big|, \Big |%
a_k \varepsilon_0 +\sum_{i \neq k} a_i \varepsilon^{\prime }_i\Big| \Big \}$%
, and taking the ${\mathbb{L}}^2$-norm in (\ref{majpo}), it follows that %
\eqref{condFM2bis} is satisfied as soon as \eqref{condcorflin} is \hfill $%
\square$

\subsection{Autoregressive Lipschitz models.}

In this section, we give an example of iterative Lipschitz model, which
fails to be irreducible, to which our results apply. For the sake of
simplicity, we do not analyze the iterative Lipschitz models in their full 
generality, as defined in Diaconis and Freedman 
\cite{DiFr} and Duflo \cite{Du}.

\medskip

For $\delta $ in $[0,1[$ and $C$ in $]0,1]$, let ${\mathcal{L}}(C,\delta )$
be the class of 1-Lipschitz functions $h$ which satisfy 
\begin{equation*}
h(0)=0\ \mbox{ and }\ |h^{\prime }(t)|\leq 1-C(1+|t|)^{-\delta }\ 
\mbox{
almost everywhere.}
\end{equation*}%
Let $(\varepsilon _{i})_{i\in {\mathbb{Z}}}$ be a sequence of iid
real-valued random variables. For $S\geq 1$, let $AR{\mathcal{L}}(C,\delta
,S)$ be the class of Markov chains on ${\mathbb{R}}$ defined by 
\begin{equation}
Y_{n}=h(Y_{n-1})+\varepsilon _{n}\ \mbox{ with }\ h\in {\mathcal{L}}%
(C,\delta )\ \mbox{ and }{\mathbb{E}}\Big (\frac{|\varepsilon
_{0}|^{S}L(\varepsilon _{0})}{L(L(\varepsilon _{0}))}\Big )<\infty \,.
\end{equation}%
(Recall that $L(x)=\log (\mathrm{e}+|x|)$). For this model, there exists a
unique invariant probability measure $\mu $ (see
Proposition 2 of Dedecker and Rio \cite{DeRi}). Moreover we have

\begin{proposition}
\label{momentofmu} Assume that $(Y_{i})_{i\in \mathbb{Z}}$ belongs to $AR{%
\mathcal{L}}(C,\delta ,S)$. Then there exists a unique invariant probability
measure that satisfies 
\begin{equation*}
\int |x|^{S-\delta }\frac{L(x)}{L(L(x))}\mu (dx)<\infty \,.
\end{equation*}
\end{proposition}

Applying Theorem 2 we derive the following result.

\begin{corollary}
\label{corarlmodel} Assume that $(Y_{i})_{i\in \mathbb{Z}}$ is a stationary
Markov chain belonging to $AR{\mathcal{L}}(C,\delta ,S)$ for some $S\geq
2+\delta $. Then, for any Lipschitz function $g$, the conclusions of Theorem
1 hold for $(g(Y_{i})-\mu (g))_{i\in \mathbb{Z}}$.
\end{corollary}

\begin{remark}
The proof of this result reveals that an application of Theorem 1 would
require the following moment condition on $\mu $: $\int |x|^{2}L(x)\mu
(dx)<\infty $ which according to the proof of Proposition \ref{momentofmu}
is satisfied provided that ${\mathbb{E}}(|\varepsilon
_{0}|^{S}L(\varepsilon _{0}))<\infty $ for some $S\geq 2+\delta $.
\end{remark}

\noindent \textbf{Proof of Proposition \ref{momentofmu}.} To prove
Proposition \ref{momentofmu}, we shall modify the proof of Proposition 2 of
Dedecker and Rio \cite{DeRi} as follows. Let $K$ be the transition kernel of
the stationary Markov chain $(Y_{i})_{i\in {\mathbb{Z}}}$ belonging to $%
ARL(C,\delta ,\eta )$. For $n>0$, we write $K^{n}g$ for the function $\int
g(y)K^{n}(x,dy)$. Let $V(x)=|x|^{S}\frac{L(x)}{L(L(x))}$. Notice that 
\begin{equation*}
KV(x)={\mathbb{E}}\big (V(Y_{n+1})|Y_{n}=x\big )={\mathbb{E}}\big (%
V(h(x)+\varepsilon _{0})\big )\leq {\mathbb{E}}\big (V(|h(x)|+|\varepsilon
_{0}|)\big )\,.
\end{equation*}%
By assumption on $h$, there exists $R_{1}\geq 1$ and some $c\in $ $]0,1/2[$
such that for every $x$, with $|x|>R_{1}$, $|h(x)|\leq |x|-c|x|^{1-\delta
}:=g(x)\leq |x|$. Therefore, using the fact that for any positive reals $a$
and $b$, $\log (\mathrm{e}+a+b)=\log (\mathrm{e}+a)+\log (1+b/(\mathrm{e}%
+a)) $, we get for any $|x|>R_{1}$ (using that for $u\geq 0$, $\log
(1+u)\leq u$), 
\begin{align}
V(|h(x)|& +|\varepsilon _{0}|)\leq (g(x)+|\varepsilon _{0}|)^{S}\frac{L(x)}{%
L(L(|x|+|\varepsilon _{0}|))}+(|x|+|\varepsilon _{0}|)^{S}\frac{%
L(\varepsilon _{0}/(1+|x|)}{L(L(|x|+|\varepsilon _{0}|))}  \notag
\label{p1propmomentmu} \\
& \leq (g(x)+|\varepsilon _{0}|)^{S}\frac{L(x)}{L(L(x))}+\frac{%
2^{S}|x|^{S-1}|\varepsilon _{0}|}{L(L(R_{1}))}+\frac{2^{S}|\varepsilon
_{0}|^{S}L(\varepsilon _{0})}{L(L(\varepsilon _{0}))}\,.
\end{align}%
To deal now with the first term in the right hand side of the above
inequality, we shall use inequality \eqref{basicine}, in the Appendix, with $%
a=g(x)$ and $b=|\varepsilon _{0}|$. We get that there exist positive
constants $c$ and $R_{2}$ such that for any $|x|>R_{2}$, 
\begin{gather*}
V(|h(x)|+|\varepsilon _{0}|)\leq (g(x))^{S}\frac{L(x)}{L(L(x))}%
+2^{S+1}|\varepsilon _{0}|(g(x))^{S-1}\frac{L(x)}{L(L(x))} \\
+2^{S}|\varepsilon _{0}|^{S}\frac{L(x)}{L(L(x))}+2^{S}|x|^{S-1}|\varepsilon
_{0}|+\frac{2^{S}|\varepsilon _{0}|^{S}L(\varepsilon _{0})}{L(L(\varepsilon
_{0}))} \\
\leq |x|^{S}\frac{L(x)}{L(L(x))}-c|x|^{S-\delta }\frac{L(x)}{L(L(x))}%
+3\times 2^{S}|\varepsilon _{0}||x|^{S-1}L(x)+2^{S}|\varepsilon _{0}|^{S}%
\frac{L(x)}{L(L(x))}+\frac{2^{S}|\varepsilon _{0}|^{S}L(\varepsilon _{0})}{%
L(L(\varepsilon _{0}))}\,.
\end{gather*}%
Taking the expectation, considering the moment assumption on $\varepsilon
_{0}$ and using the fact that $\delta \in \lbrack 0,1[$ and $S\geq 1$, it
follows that there exist positive constants $d$ and $R$ such that for any $%
|x|>R$, 
\begin{equation*}
KV(x)\leq V(x)-d\,|x|^{S-\delta }\frac{L(x)}{L(L(x))}\,.
\end{equation*}%
So overall it follows that there exists a positive constant $b$ such that 
\begin{equation}
KV(x)\leq V(x)-d\,|x|^{S-\delta }\frac{L(x)}{L(L(x))}+b\mathbf{1}%
_{[-R,R]}(x)\,.  \label{p2propmomentmu}
\end{equation}%
This inequality allows to use the arguments given at the end of the proof of
Proposition 2 in Dedecker and Rio \cite{DeRi}. Indeed, iterating $n$ times the inequality %
\eqref{p2propmomentmu}, we get 
\begin{equation*}
\frac{d}{n}\sum_{k=1}^{n}\int |y|^{S-\delta }\frac{L(y)}{L(L(y))}%
K^{k}(x,dy)\leq \frac{1}{n}KV(x)+\frac{b}{n}\sum_{k=1}^{n}K^{k}([-R,R])(x)\,,
\end{equation*}%
and letting $n$ tend to infinity, it follows that 
\begin{equation*}
d\int |x|^{S-\delta }\frac{L(x)}{L(L(x))}\mu (dx)\leq b\mu ([-R,R])<\infty
\,.
\end{equation*}%
$\square $

\medskip

\noindent \textbf{Proof of Corollary \ref{corarlmodel}.} Let $X_i = g(Y_i) -
\mu (g)$. According to Proposition \ref{momentofmu}, ${\mathbb{E}} \Big (
\frac{X_0^2 L( X_0)}{L(L( X_0))} \Big ) < \infty$. Hence Corollary \ref%
{corarlmodel} will follow from Theorem \ref{thmlilmoment2delta} if we can
prove that the condition \eqref{condFM} is satisfied which will clearly hold
if 
\begin{equation}  \label{p1corarlmodel}
\sum_{n >0} n^{-1}\vert K^n g(x) - \mu (g) |^2 \mu (dx) < \infty \, .
\end{equation}
According to the inequality (5.7) in Dedecker and Rio \cite{DeRi}, there
exists a positive constant $A$ such that 
\begin{multline}  \label{p2corarlmodel}
\vert K^n g(x) - \mu (g) | \leq A n^{1-S/\delta} \int|x-y|\mu(dy)+ A
(1-B_n(x))^n\int|x-y|\mu(dy) \\
+A n^{1-(S-1)/\delta} \int|x-y| |y|^{S-\delta - 1}\mu(dy) \, ,
\end{multline}
where $B_n(x)=C[4(1+|x|+(n-1){{\mathbb{E}}}|\varepsilon_0|)]^{-\delta}$.
Noticing that $\sum_{n>0} n^{1-2(S-1)/\delta} < \infty$ as soon as $S >1+
\delta$ and that according to Proposition \ref{momentofmu}, $x^2$ is $\mu-$%
integrable as soon as $S \geq 2 + \delta$, we infer from %
\eqref{p2corarlmodel} that \eqref{p1corarlmodel} will be satisfied if we can
prove that 
\begin{equation}  \label{p3corarlmodel}
\sum_{n >0} n^{-1} \int x^2 (1-B_n(x))^{2n}\mu(dx) < \infty \, .
\end{equation}
Notice that that $(1-B_n(x))^{2n}\leq \exp(-2nB_n(x))$. If $|x| \leq 1 $, 
\begin{equation*}
\exp(-2nB_n(x)) \leq \exp ( -2 Cn (8 + 4 n )^{-\delta} ) \, ,
\end{equation*}
implying that 
\begin{equation*}
\int_{-1}^1 x^2 \sum_{n >0} n^{-1} \exp(-2nB_n(x))\mu(dx) < \infty \, .
\end{equation*}
Now if $|x| >1$, 
\begin{align*}
\sum_{n \geq 2} n^{-1} \exp(-2nB_n(x)) & \leq \int_1^{\infty} u^{-1}\exp %
\big ( - 2 C u [4(1+|x|+u {{\mathbb{E}}}|\varepsilon_0|)]^{-\delta} \big ) du
\\
& \leq \int_1^{\infty} u^{-1}\exp \big ( - 2 C u |x|^{-\delta }[ 8 + u
|x|^{-\delta} {{\mathbb{E}}}|\varepsilon_0|)]^{-\delta} \big ) du \\
& \leq \int_1^{\infty} z^{-1}\exp \big ( - 2 C z[ 8 + z {{\mathbb{E}}}%
|\varepsilon_0|)]^{-\delta} \big ) dz \, .
\end{align*}
Hence there exists a positive constant $M$ such that 
\begin{equation}  \label{p4corarlmodel}
\int_{|x| >1} x^2 \sum_{n >0} n^{-1} \exp(-2nB_n(x))\mu(dx) \leq M \int_{|x|
>1} x^2 \mu (dx) \, ,
\end{equation}
which according to Proposition \ref{momentofmu} is finite as soon as $S \geq
2 + \delta$. All the above computations then show that \eqref{p3corarlmodel}
(and then \eqref{p1corarlmodel}) holds provided that $S \geq 2 + \delta$. $%
\square$

\subsection{Application to weakly dependent sequences}

Theorems \ref{thmlilmoment2} and \ref{thmlilmoment2delta} can be
successively applied to large classes of weakly dependent sequences. In this
section, we give an application to $\alpha$-dependent sequences. With this
aim, we first need some definitions.

\begin{definition}
\label{defalphaweak} For a sequence $\mathbf{Y}=(Y_{i})_{i\in {\mathbb{Z}}}$%
, where $Y_{i}=Y_{0}\circ \theta ^{i}$ and $Y_{0}$ is an $\mathcal{F}_{0}$%
-measurable and real-valued random variable, let for any $k\in {\mathbb{N}}$%
, 
\begin{equation*}
\alpha _{{\mathbf{Y}}}(k)=\sup_{u\in {\mathbb{R}}}\big \|\mathbb{E}(\mathbf{1%
}_{Y_{k}\leq u}|{\mathcal{F}}_{0})-\mathbb{E}(\mathbf{1}_{Y_{k}\leq u})\big
\|_{1}\,.
\end{equation*}
\end{definition}

\begin{remark}
\label{defalphastrong} Recall that the strong mixing coefficient of
Rosenblatt \cite{Ros56} between two $\sigma$-algebras ${\mathcal{F}}$ and ${%
\mathcal{G}}$ is defined by 
\begin{equation*}
\alpha({\mathcal{F}}, {\mathcal{G}})= \sup_{A \in {\mathcal{F}}, B \in {%
\mathcal{G}}}|{\mathbb{P}}(A \cap B)-{\mathbb{P}}(A){\mathbb{P}}(B)| \, .
\end{equation*}
For a strictly stationary sequence $(Y_i)_{i \in {\mathbb{Z}}}$ of real
valued random variables, and the $\sigma$-algebra ${\mathcal{F}}_0=\sigma
(Y_i, i \leq 0)$, define then 
\begin{equation}  \label{defalpharosen}
\alpha(0) = 1 \text{ and } \alpha(k)= 2\alpha({\mathcal{F}}_0,\sigma(Y_k)) 
\text{ for $k>0$} \, .
\end{equation}
Between the two above coefficients, the following relation holds: for any
positive $k$, $\alpha_{ \mathbf{Y}}(k) \leq \alpha (k)$. In addition, the $%
\alpha$-dependent coefficient as defined in Definition \ref{defalphaweak}
may be computed for instance for many Markov chains associated to dynamical
systems that fail to be strongly mixing in the sense of Rosenblatt \cite%
{Ros56}.
\end{remark}

\begin{definition}
\label{defquant} A quantile function $Q$ is a function from $]0,1]$ to ${%
\mathbb{R}}_+$, which is left-continuous and non increasing. For any
nonnegative random variable $Z$, we define the quantile function $Q_Z $ of $%
Z $ by $Q_Z (u) = \inf \{ t \geq 0 : {\mathbb{P}} (|Z| >t ) \leq u \} $.
\end{definition}

\begin{definition}
\label{def1} Let $\mu$ be the probability distribution of a random variable $%
X$. If $Q$ is an integrable quantile function (see Definition \ref{defquant}%
), let $\Mon( Q, \mu)$ be the set of functions $g$ which are monotonic on
some open interval of ${\mathbb{R}}$ and null elsewhere and such that $%
Q_{|g(X)|} \leq Q$. Let $\Mon^c( Q, \mu)$ be the closure in ${\mathbb{L}}%
^1(\mu)$ of the set of functions which can be written as $\sum_{\ell=1}^{L}
a_\ell f_\ell$, where $\sum_{\ell=1}^{L} |a_\ell| \leq 1$ and $f_\ell$
belongs to $\Mon( Q, \mu)$.
\end{definition}

Applying Theorem \ref{thmlilmoment2delta}, we get

\begin{corollary}
\label{coralpha} Let $Y_0$ be a real-valued random variable with law $%
P_{Y_0} $, and $Y_i=Y_0 \circ \theta^i$. Let $X_i = f(Y_i) - \mathbb{E} (
f(Y_i))$ where $f$ belongs to $\Mon^c( Q, P_{Y_0})$ with $Q^2 L(Q)/ L(L(Q)))$
integrable. Assume in addition that 
\begin{equation}
\sum_{k\geq 3} \frac{1}{ k (\log \log k )}\int_0^{\alpha_{\mathbf{Y}} (k)}
Q^2 (u) du <\infty \, .  \label{condalphaQ}
\end{equation}%
Then (\ref{condFM}) is satisfied and consequently, the conclusions of
Theorem \ref{thmlilmoment2} hold for $(X_k)_{k \in {\mathbb{Z}}}$.
\end{corollary}

To prove that (\ref{condalphaQ}) implies (\ref{condFM}), it suffices to
notice that 
\begin{equation*}
\Vert {\mathbb{E}}_0(X_k) \Vert_2^2 = {\mathbb{E}} ( X_k {\mathbb{E}}_0 (
X_k) ) \leq \int_0^{\alpha_{\mathbf{Y}} (k)} Q^2 (u) du
\end{equation*}
(see the proof of (4.17) in Merlev\`ede and Rio \cite{MeRi} for the last
inequality).

\medskip

The definition \ref{def1} describes spaces similar to weak ${\mathbb{L}}^p$
where we require a monotonicity condition plus a uniform bound on the tails
of the functions. Let us introduce in the same spirit ${\mathbb{L}}^p$-like
spaces.

\begin{definition}
\label{def2}If $\mu$ is a probability measure on ${\mathbb{R}}$ and $p \in
[2,\infty)$, $M \in (0, \infty)$, let $\Mon_p(M,\mu)$ denote the set of
functions $f:{\mathbb{R}} \rightarrow {\mathbb{R}} $ which are monotonic on
some interval and null elsewhere and such that $\mu(|f|^p)\leq M^p$. Let $%
\Mon^c_p(M,\mu)$ be the closure in ${\mathbb{L}}^p(\mu)$ of the set of
functions which can be written as $\sum_{\ell=1}^L a_\ell f_\ell$, where $%
\sum_{\ell=1}^L|a_\ell| \leq 1$ and $f_\ell\in \Mon_p(M, \mu)$.
\end{definition}

Let $X_i = f(Y_i) - \mathbb{E} ( f(Y_i))$, where $f$ belongs to $\Mon^c_{2 +
\delta}(M,P_{Y_0})$ for some $\delta >0$. If 
\begin{equation*}
\sum_{k\geq 3}\frac{(\alpha_{\mathbf{Y}} (k))^{\delta/(2+\delta)}}{k (\log
\log k)}<\infty \, ,  \label{condalpha}
\end{equation*}%
then Corollary \ref{coralpha} applies.

\medskip

\noindent \textbf{Application to functions of Markov chains associated to
intermittent maps.}

\medskip

For $\gamma$ in $]0, 1[$, we consider the intermittent map $T_\gamma$ from $%
[0, 1]$ to $[0, 1]$, which is a modification of the Pomeau-Manneville map 
\cite{PoMa}: 
\begin{equation*}
T_\gamma(x)= 
\begin{cases}
x(1+ 2^\gamma x^\gamma) \quad \text{ if $x \in [0, 1/2[$} \\ 
2x-1 \quad \quad \quad \ \ \text{if $x \in [1/2, 1]$} \, .%
\end{cases}%
\end{equation*}
Recall that $T_{\gamma}$ is ergodic and that there exists a unique $T_\gamma$%
-invariant probability measure $\nu_\gamma$ on $[0, 1]$, which is absolutely
continuous with respect to the Lebesgue measure. We denote by $L_\gamma$ the
Perron-Frobenius operator of $T_\gamma$ with respect to $\nu_\gamma$. Recall
that for any bounded measurable functions $f$ and $g$, 
\begin{equation*}
\nu_\gamma(f \cdot g\circ T_\gamma)=\nu_\gamma(L_\gamma(f) g) \, .
\end{equation*}
Let $(Y_i)_{i \geq 0}$ be a Markov chain with transition Kernel $L_\gamma$
and invariant measure $\nu_\gamma$.

\begin{corollary}
\label{ASmapB} Let $\gamma \in (0,1)$ and $(Y_i)_{i \geq 1}$ be a stationary
Markov chain with transition kernel $L_{\gamma}$ and invariant measure $%
\nu_\gamma$. Let $Q$ be a quantile function such that 
\begin{equation}  \label{lilcondIT}
\int_0^1 \frac{L (Q(u))}{L ( L (Q(u))} Q^2 (u) du <\infty\,.
\end{equation}
Let $X_i=f(Y_i) - \nu_{\gamma}(f)$ where $f$ belongs to $\Mon^c( Q,
\nu_{\gamma})$. Then (\ref{condFM}) is satisfied and the conclusions of
Theorem \ref{thmlilmoment2} hold for $(X_k)_{k \in {\mathbb{Z}}}$.
\end{corollary}

\begin{remark}
Notice that, by standard arguments on quantile functions, \eqref{lilcondIT} is equivalent to the following condition: 
\begin{equation*}
\int_{0}^{\infty }\frac{xL(x)}{L(L(x))}Q^{-1}(x)dx<\infty \,,
\end{equation*}%
where $Q^{-1}$ is the generalized inverse of $Q$.
\end{remark}

\noindent \textbf{Proof.} To prove this corollary, it suffices to see that (%
\ref{lilcondIT}) implies (\ref{condalphaQ}). For this purpose, we first
notice that \eqref{condalphaQ} can be rewritten in the following equivalent
way (see Rio \cite{Rio}): 
\begin{equation*}
\int_{0}^{1}\frac{L(\alpha _{\mathbf{Y}}^{-1}(u))}{L(L(\alpha _{\mathbf{Y}%
}^{-1}(u)))}Q^{2}(u)du<\infty \,,
\end{equation*}%
where $\alpha _{\mathbf{Y}}^{-1}(x)=\min \{q\in {\mathbb{N}}\,:\,\alpha _{%
\mathbf{Y}}(q)\leq x\}$. Now, according to Proposition 1.17 in Dedecker 
\textit{et al.} \cite{DGM}, there exists a positive constant $C$ such that $%
\alpha _{\mathbf{Y}}^{-1}(u)\leq Cu^{-\gamma /(1-\gamma )}$. Therefore, for
any $\eta \in ]0,1/2[$, there exists a constant $c$ depending on $\gamma $, $%
C$ and $\eta $ such that 
\begin{gather*}
\int_{0}^{1}\frac{L(\alpha _{\mathbf{Y}}^{-1}(u))}{L(L(\alpha _{\mathbf{Y}%
}^{-1}(u)))}Q^{2}(u)du\leq c\int_{0}^{1}\frac{L(u^{-\eta })}{L(L(u^{-\eta }))%
}Q^{2}(u)du \\
\leq c\int_{0}^{1}\frac{L(Q(u))}{L(L(Q(u))}Q^{2}(u)du+c\int_{0}^{1}\frac{%
L(u^{-\eta })}{L(L(u^{-\eta }))}u^{-2\eta }du\,,
\end{gather*}%
which is finite under \eqref{condalphaQ}. $\square $

\medskip

In particular, if $f$ is with bounded variation and $\gamma<1$, we infer
from Corollary \ref{ASmapB} that the conclusions of Theorem \ref%
{thmlilmoment2} hold for $(X_k)_{k \in {\mathbb{Z}}}$. Note also that %
\eqref{lilcondIT} is satisfied if $Q$ is such that $Q(u)\leq C u^{-1/2}(
\log (u^{-1} ))^{-1}( \log \log (u^{-1}))^{-b}$ for $u$ small enough and $%
b>1/2$. Therefore, since the density $h_{\nu_{\gamma}}$ of $\nu_{\gamma}$ is
such that $h_{\nu_{\gamma}}(x) \leq C x^{-\gamma}$ on $(0, 1]$, one can
easily prove that if $f$ is positive and non increasing on $(0, 1)$, with 
\begin{equation*}
f(x) \leq \frac{C}{x^{(1 - \gamma)/2}| \log (x) |(\log |\log (x)| )^{b}}
\quad \text{near 0 for some $b>1/2$},
\end{equation*}
then \eqref{lilcondIT} holds.

\subsection{Application to a class of Markov chains}

We shall point out next some consequences of Theorem \ref{thm4} in terms of
stationary Markov chains characteristics. Let $T$ be a regular transition
probability on the measurable space $({\mathbb{S}},\mathcal{S})$, leaving
invariant a probability $\mu $ on $({\mathbb{S}},\mathcal{S)}$. We also
denote by $T$ the Markov operator induced on ${\mathbb{L}}^{2}(\mu )$ via%
\begin{equation*}
Tg(x)=\int_{{\mathbb{S}}}g(y)T(x,dy)
\end{equation*}%
and we assume that it is ergodic (i.e. $Tf=f\ $ $\mu $ a.e. for $f\geq 0$
implies $f$ is constant $\mu $ a.e.). Let $(\xi _{n})_{n\in \mathbb{Z}}$ be
the (stationary) canonical Markov chain with state space $({\mathbb{S}},%
\mathcal{S})$ associated with $T$, defined on the canonical space $(\Omega ,{%
\mathcal{A}},{\mathbb{P}})=({\mathbb{S}}^{\mathbb{Z}},\mathcal{S}^{\otimes 
\mathbb{Z}},{\mathbb{P}})$ (the law of $\xi _{0}$ under ${\mathbb{P}}$ is $%
\mu $). Denote by $\theta $ the shift on $\Omega .$ Denote by $\mathcal{F}%
_{n}=\sigma (\xi _{j};j\leq n).$

\begin{corollary}
\label{MarkovChains} Let $(\xi _{n})_{n\in \mathbb{Z}}$ be a stationary
Markov chain. Let $h\in {\mathbb{L}}^{2}({\mathbb{S}},\mu )$ centered and
define $X_{k}=h(\xi _{k})$. Let $t\in (0,2\pi )\backslash \{\pi \}$ be such
that $\mathrm{e}^{-it}$ is not in the spectrum of $T$, and $\mathrm{e}%
^{-2it} $ is not an eigenvalue of $T$. Then the conclusions of Theorem \ref%
{thm4} hold for $(X_{k})_{k\in {\mathbb{Z}}}$.
\end{corollary}

\begin{remark}
In Corollary \ref{MarkovChains} we do not assume the regularity condition $%
\Vert T^{n} h \Vert _{2,\mu }\rightarrow 0$. The spectral density might not
exist.
\end{remark}

\noindent \textbf{Proof of Corollary \ref{MarkovChains}.} By assumption,
there exists $g\in {\mathbb{L}}^{2}(\SS,\mu )$ such that $h=g-\mathrm{e}%
^{it}Tg$. Therefore, 
\begin{equation*}
\mathbb{E}_{0}(S_{n}(t))=\sum_{k=1}^{n}\mathbb{E}_{0}(\mathrm{e}^{itk}g(\xi
_{k})-\mathrm{e}^{it(k+1)}g(\xi _{k+1}))=\mathrm{e}^{it}\mathbb{E}_{0}(g(\xi
_{1}))-\mathrm{e}^{it(n+1)}\mathbb{E}_{0}(g(\xi _{n+1}))\,,
\end{equation*}%
showing that condition (\ref{poisson}) is satisfied. Hence, for $t\in
(0,2\pi )\backslash \{\pi \}$ such that $\mathrm{e}^{-2it}$ is not an
eigenvalue of $\theta $ (the shift on $\Omega $), the proposition will
follow from Theorem \ref{thm4}. To end the proof we notice that if $\mathrm{e%
}^{-2it}$ is an eigenvalue for $\theta $, it is an eigenvalue for $T$, see
e.g. Proposition 2.3 of Cuny \cite{Cuny0} (notice that the proof there
extends easily to ${\mathbb{L}}^{2}$)\hfill $\square $

\medskip

We give now a consequence of Corollary \ref{MarkovChains} for reversible
Markov chains. This follows from the fact that the spectrum of $T$ is real
and lies in $[-1,1]$.

\begin{corollary}
\label{corRevMC} Let $(\xi _{n})_{n\in \mathbb{Z}}$ be a stationary and
reversible Markov chain $(T=T^{\ast })$. Let $h\in {\mathbb{L}}^{2}({\mathbb{%
S}},\mu )$ centered and define $X_{k}=h(\xi _{k})$. Then the conclusion of
Theorem \ref{thm4} holds for every $t\in (0,2\pi )\backslash \{\pi \}.$
\end{corollary}

The (independent) Metropolis Hastings Algorithm leads to a Markov chain with
transition function 
\begin{equation*}
T(x,A)=p(x)\delta _{x}(A)+(1-p(x))\nu (A)\,,
\end{equation*}%
where $\delta _{x}$ denotes the Dirac measure at point $x$, $\nu $ is a
probability measure on ${\mathcal{S}}$ and $p:{\mathbb{S}}\rightarrow
\lbrack 0,1]$ is a measurable function for which $\theta =\int_{{\mathbb{S}}}%
\frac{1}{1-p(x)}\nu (dx)<\infty $. Then there is a unique invariant
distribution 
\begin{equation*}
\mu (dx)=\frac{1}{\theta (1-p(x))}\nu \,(dx)
\end{equation*}%
and the associated stationary Markov chain $(\xi _{i})_{i}$ is reversible
and ergodic. Hence Corollary \ref{corRevMC} applies to this example.

\section{Appendix}

\subsection{Facts from ergodic theory}

We first recall the following consequence of the Dunford-Schwartz ergodic
theorem, see sections VIII.5 and VIII.6 of \cite{DS}.

\begin{proposition}
Let $(\Omega ,{\mathcal{A}},{\mathbb{P}})$ be a probability space and $%
\theta $ be a measure-preserving transformation on $\Omega $. Let $s\in 
\mathbb{R}$. For every $X\in {\mathbb{L}}^{1}(\Omega ,{\mathcal{A}},{\mathbb{%
P}})$, there exists $\pi _{s}(X)\in {\mathbb{L}}^{1}(\Omega ,{\mathcal{A}},{%
\mathbb{P}})$ such that 
\begin{equation}
\frac{1}{n}\sum_{k=0}^{n-1}\mathrm{e}^{iks}X\circ \theta ^{k}\underset{%
n\rightarrow \infty }{\longrightarrow }\pi _{s}(X)\qquad 
\mbox{${\mathbb
P}$-a.s.}  \label{DS}
\end{equation}%
and in ${\mathbb{L}}^{1}(\Omega ,{\mathcal{A}},{\mathbb{P}})$. Moreover, $%
\pi _{s}(X)\circ \theta =\mathrm{e}^{-is}\pi _{s}(X)$ ${\mathbb{P}}$-a.s.
\end{proposition}

\begin{remark}
It follows from the Wiener-Wintner theorem that the set of measure $1$ in %
\eqref{DS} may be chosen independently of $s$, but we shall not need that
refinement.
\end{remark}

\noindent \textbf{Proof.} Define an operator $V_{s}$ on ${\mathbb{L}}%
^{1}(\Omega ,{\mathcal{A}},{\mathbb{P}})$, by $V_{s}(X)=\mathrm{e}%
^{is}X\circ \theta $. Then, $V_{s}$ is a contraction of ${\mathbb{L}}^{1}$
which also contracts the ${\mathbb{L}}^{\infty }$ norm. Hence we may apply 
\cite[Theorem 6 p. 675]{DS}, to obtain the almost sure$\,$convergence. The ${%
\mathbb{L}}^{1}$ convergence follows from \cite[Corollary 5 p. 664]{DS} (see
also the proof of the next lemma). \hfill $\square $

\medskip

We also give the following lemma, that should be well-known. We give a proof
for completeness.

\begin{lemma}
\label{LLNFourier} Let $(\Omega ,{\mathcal{A}},{\mathbb{P}})$ be a
probability space and $\theta $ be a measure-preserving transformation on $%
\Omega $. Let $t_{0}\in \mathbb{R}$ be fixed. If there is no non trivial $%
Y\in {\mathbb{L}}^{2}(\Omega ,{\mathcal{A}},{\mathbb{P}})$, such that $%
Y\circ \theta =\mathrm{e}^{-it_{0}}Y$ \mbox{${\mathbb P}$-a.s.}, then, for
every $X\in {\mathbb{L}}^{1}(\Omega ,{\mathcal{F}},{\mathbb{P}})$ $\pi
_{t_{0}}(X)=0$ \mbox{${\mathbb P}$-a.s.}\thinspace Furthermore, when ${%
\mathbb{L}}^{2}(\Omega ,{\mathcal{A}},{\mathbb{P}})$ is separable, there
exists a countable (at most) set $\mathbf{S}\subset \mathbb{R}$ such that
for every $t\in \mathbb{R}\backslash \mathbf{S}$ and every $X\in {\mathbb{L}}%
^{1}(\Omega ,{\mathcal{A}},{\mathbb{P}})$, $\pi _{t}(X)=0$ 
\mbox{${\mathbb
P}$-a.s.}
\end{lemma}

\noindent \textbf{Proof.} Define $V_{t_{0}}$ as above. Since $\sup_{n\geq 1}%
\frac{1}{n}\Vert \sum_{k=0}^{n-1}V_{t_{0}}^{k}\Vert _{{\mathbb{L}}%
^{1}\rightarrow {\mathbb{L}}^{1}}<\infty $, by the Banach-principle, the set 
${\mathcal{Y}}:=\{X\in {\mathbb{L}}^{1}(\Omega ,{\mathcal{A}},{\mathbb{P}}%
)~:~\Vert \frac{1}{n}\sum_{k=0}^{n-1}\mathrm{e}^{ikt_{0}}X\circ \theta
^{k}\Vert _{1}\rightarrow 0\}$ is closed in ${\mathbb{L}}^{1}$. Now, by von
Neumann's mean ergodic theorem 
\begin{equation*}
{\mathbb{L}}^{2}(\Omega ,{\mathcal{A}},{\mathbb{P}})=\overline{(I-V_{t_{0}}){%
\mathbb{L}}^{2}(\Omega ,{\mathcal{A}},{\mathbb{P}})}\oplus \mathrm{Fix}%
\,V_{t_{0}},
\end{equation*}%
where $\mathrm{Fix}\,V_{t_{0}}$ stands for the fixed points of $V_{t_{0}}$
in ${\mathbb{L}}^{2}(\Omega ,{\mathcal{A}},{\mathbb{P}})$ and the closure is
in norm $\Vert \cdot \Vert _{2}$. By assumption, $V_{t_{0}}$ has no non
trivial fixed point. Obviously, ${\mathcal{Y}}$ contains $(I-V_{t_{0}}){%
\mathbb{L}}^{2}(\Omega ,{\mathcal{A}},{\mathbb{P}})$, hence ${\mathcal{Y}}={%
\mathbb{L}}^{1}(\Omega ,{\mathcal{A}},{\mathbb{P}})$. \newline
Assume now that ${\mathbb{L}}^{2}(\Omega ,{\mathcal{A}},{\mathbb{P}})$ is
separable. Define an operator on ${\mathbb{L}}^{2}(\Omega ,{\mathcal{A}},{%
\mathbb{P}})$ by $UX=X\circ \theta $. It is well known that the eigenspaces
of $U$ corresponding to different eigenvalues are orthogonal. By
separability there are at most countably many eigenvalues for $U$, hence the
result. \hfill $\square $

\subsection{Technical approximation results}

\begin{lemma}
\label{max2} Assume that $X_0$ is almost surely bounded by $M$. For any
integer $s\geq 1$%
\begin{equation*}
(2 \pi)^{-1} \int_0^{2\pi}\mathbb{E}(\max_{1\leq \ell\leq m}|S_{\ell}(t
)-M_{\ell}(t )|^{2}) dt \leq 12\big ( m \Vert \mathbb{E}_{-s}(X_{0})
\Vert^{2}+s^{2}M^{2} \big ) \, ,
\end{equation*}
where $M_{n}(t)=\sum_{k=1}^{n}D_{k}(t)$ and $D_{k}(t)$ is defined by %
\eqref{defDlalpha}.
\end{lemma}

The proof of this lemma follows from the following result (which is of
independent interest) by selecting $a_{j}=1$ for $0\leq j\leq s-1$ and $%
a_{j}=0$ for any $j\geq s$.

\begin{lemma}
\label{max1}Assume that $X_0$ is almost surely bounded by $M$. Let $(a_{n})$
be a sequence of positive numbers nonincreasing smaller than $1$ with $%
\sum_{j=1}^{\infty }a_{j}<\infty $ and $a_{0}=1$. Then 
\begin{equation*}
(2 \pi)^{-1}\int_0^{2\pi}\mathbb{E}(\max_{1\leq \ell\leq m}|S_{\ell}(t
)-M_{\ell}(t )|^{2}) dt \leq 12 \Big ( m\sum_{j=1}^{\infty }(a_{j-1}-a_{j})
\Vert \mathbb{E}_{-j}(X_{0}) \Vert^{2}+M^{2} \big ( \sum_{j=0}^{\infty
}a_{j} \big )^{2} \Big ) \, .
\end{equation*}
\end{lemma}

\noindent \textbf{Proof of Lemma \ref{max1}.}

\medskip

\noindent \textit{Step 1: Martingale decomposition.}

We start with a traditional martingale decomposition (see for instance
Section 4.1 in Merlev\`{e}de, Peligrad and Utev \cite{MPU}). Let $t\in
\lbrack 0,2\pi )$ and $X_{k}(t)=\mathrm{e}^{ikt}X_{k}$. 
\begin{equation*}
\theta _{k}(t)=X_{k}(t)+\sum_{j=1}^{\infty }a_{j}\mathbb{E}_{k}(X_{k+j}(t));%
\text{ \ \ \ \ }\theta _{k}^{\prime }(t)=\sum_{j=1}^{\infty }a_{j\ }\mathbb{E%
}_{k}(X_{k+j}(t))
\end{equation*}%
and 
\begin{equation*}
\mathbb{E}_{k}(\theta _{k+1}(t))-\theta _{k}(t)=-X_{k}(t)+\sum_{j=1}^{\infty
}(a_{j-1}-a_{j})\mathbb{E}_{k}(X_{k+j}(t))\,.
\end{equation*}%
Finally, denote by 
\begin{equation*}
D_{k+1}^{\prime }(t)=\theta _{k+1}(t)-\mathbb{E}_{k}(\theta
_{k+1})(t)=\sum_{j=0}^{\infty }a_{j}P_{k+1}(X_{k+j+1}(t))\ \ ;\ \
M_{n}^{\prime }(t)=\sum_{k=1}^{n}D_{k}^{\prime }(t)\,.
\end{equation*}%
Then, $(D_{k}^{\prime }(t))_{k\in \mathbb{Z}}$ is a sequence of martingale
differences with respect to the stationary filtration $({\mathcal{F}}%
_{j})_{j\in \mathbb{Z}}.$ Note 
\begin{equation*}
X_{k}(t)=D_{k+1}^{\prime }(t)+\theta _{k}(t)-\theta
_{k+1}(t)+\sum_{j=1}^{\infty }(a_{j-1}-a_{j})\mathbb{E}_{k}(X_{k+j}(t))\,.
\end{equation*}%
Taking into account the definition of $\theta _{k}^{^{\prime }}(t)$ we can
also write%
\begin{equation*}
X_{k}(t)=D_{k}^{\prime }(t)+\theta _{k-1}^{\prime }(t)-\theta _{k}^{\prime
}(t)+\sum_{j=1}^{\infty }(a_{j-1}-a_{j})\mathbb{E}_{k-1}(X_{k+j-1}(t))\,.
\end{equation*}%
It follows that for almost all $t\in \lbrack 0,2\pi )$, 
\begin{eqnarray}
S_{\ell }(t)-M_{\ell }(t) &=&\sum_{k=1}^{\ell }\sum_{j=0}^{\infty
}(a_{j}-1)P_{k}(X_{k+j}(t))+\sum_{k=0}^{\ell -1}\sum_{j=1}^{\infty
}(a_{j-1}-a_{j})\mathbb{E}_{k}(X_{k+j}(t))+\theta _{0}^{\prime }(t)-\theta
_{\ell }^{\prime }(t)  \notag \\
&=&I+I\!I+\theta _{0}^{\prime }(t)-\theta _{m}^{\prime }(t)\,.
\label{defrest}
\end{eqnarray}%
\textit{Step 2: The estimation of } $\int_{0}^{2\pi }\mathbb{E}\max_{1\leq
\ell \leq m}|S_{\ell }(t)-M_{\ell }(t)|^{2}dt.$

We shall estimate separately this maximum for all the terms in the
decomposition \eqref{defrest}.

By the Doob-Kolmogorov martingale maximal inequality, stationarity, Fubini
theorem and orthogonality of $\mathrm{e}^{itk}$ the first term will be
(remind $a_{0}=1)$ dominated by%
\begin{eqnarray*}
(2\pi )^{-1}\int_{0}^{2\pi }\mathbb{E}\max_{1\leq \ell \leq m}\mathbf{|}%
I|^{2}dt &\leq &(2\pi )^{-1}m\int_{0}^{2\pi }\mathbb{E}|D_{0}^{\prime
}(t)-D_{0}(t)|^{2}dt \\
&=&(2\pi )^{-1}m\int_{0}^{2\pi }{\mathbb{E}}\big |\sum_{j=0}^{\infty
}(a_{j}-1)P_{0}(X_{j}(t))\big |^{2}dt \\
&=&m\sum_{j=1}^{\infty }(a_{j}-1)^{2}\Vert \mathbb{P}_{-j}X_{0}|\Vert ^{2}\,.
\end{eqnarray*}%
By simple computations, 
\begin{gather*}
\sum_{j=1}^{\infty }(a_{j}-1)^{2}\Vert \mathbb{P}_{-j}X_{0}|\Vert
^{2}=\sum_{j=1}^{\infty }(a_{j}-1)^{2}\big (\Vert \mathbb{E}%
_{-j}(X_{0})\Vert ^{2}-\Vert \mathbb{E}_{-j-1}(X_{0})\Vert ^{2}\big ) \\
=(a_{1}-1)^{2}\Vert \mathbb{E}_{-1}(X_{0})\Vert ^{2}+\sum_{j=2}^{\infty
}[(a_{j}-1)^{2}-(a_{j-1}-1)^{2}]\Vert \mathbb{E}_{-j}(X_{0})\Vert ^{2} \\
=\sum_{j=1}^{\infty }[(a_{j}-1)^{2}-(a_{j-1}-1)^{2}]\Vert \mathbb{E}%
_{-j}(X_{0})\Vert ^{2}\leq 2\sum_{j=1}^{\infty }(a_{j-1}-a_{j})\Vert \mathbb{%
E}_{-j}(X_{0})\Vert ^{2}\,.
\end{gather*}%
So%
\begin{equation}
(2\pi )^{-1}\int_{0}^{2\pi }\mathbb{E}\max_{1\leq \ell \leq m}|I|^{2}dt\leq
2m\sum_{j=1}^{\infty }(a_{j-1}-a_{j})\Vert \mathbb{E}_{-j}(X_{0})\Vert
^{2}\,.  \label{rel1}
\end{equation}%
By Cauchy Schwarz inequality the second term is estimated as follows.%
\begin{equation*}
\int_{0}^{2\pi }\mathbb{E}\max_{1\leq \ell \leq m}|I\!I|^{2}dt\leq
\sum_{\ell =1}^{\infty }(a_{\ell -1}-a_{\ell })\sum_{j=1}^{\infty
}(a_{j-1}-a_{j})\int_{0}^{2\pi }\mathbb{E}\max_{1\leq \ell \leq m}\big |%
\sum_{k=0}^{\ell -1}\mathrm{e}^{i(k+j)t}\mathbb{E}_{k}(X_{k+j})\big |%
^{2}dt\,.
\end{equation*}%
By Hunt and Young maximal inequality \cite{HY}, 
\begin{align}
(2\pi )^{-1}\int_{0}^{2\pi }\mathbb{E}\max_{1\leq \ell \leq m}|I\!I|^{2}dt&
\leq \sum_{j=1}^{\infty }(a_{j-1}-a_{j})\sum_{k=0}^{m-1}\Vert {\mathbb{E}}%
_{0}(X_{j})\Vert ^{2}  \notag  \label{rel2} \\
& =m\sum_{j=1}^{\infty }(a_{j-1}-a_{j})\Vert {\mathbb{E}}_{0}(X_{j})\Vert
^{2}\,.
\end{align}%
The last terms are estimated in a trivial way as follows: 
\begin{equation}
(2\pi )^{-1}\int_{0}^{2\pi }\mathbb{E}\max_{1\leq \ell \leq m}\big |\theta
_{0}^{\prime }(t)-\theta _{k}^{\prime }(t)\big |^{2}dt=4M^{2}\big (%
\sum_{j=0}^{\infty }a_{j}\big )^{2}\,.  \label{rel3}
\end{equation}%
Gathering \eqref{rel1}, \eqref{rel2} and \eqref{rel3}, the lemma follows.
\hfill $\square $

\subsection{An algebraic inequality}

\begin{lemma}
\label{alg}For any positive reals $a$ and $b$ and any real $S\geq 1$, 
\begin{equation}
(a+b)^{S}\leq 2^{S}b^{S}+a^{S}(1+2^{S+1}b/a)\,.  \label{basicine}
\end{equation}
\end{lemma}

\noindent \textbf{Proof.} To prove the above inequality we first notice that
if $a\leq b$, the inequality is trivial. Let then assume that $b<a$. The
Newton binomial formula gives 
\begin{gather*}
(a+b)^{S}\le a^S(1+b/a)^{[S]+1} \le a^S \big ( 1+b/a\sum_{k=1}^{[S]+1} C_{[S]+1}^k (b/a)^{k-1} \big )
\le a^S (1+2^{S+1}b/a)\, .
\end{gather*}
\hfill $\square $

\end{document}